\documentclass[11pt,leqno,a4]{amsart}

\usepackage[all, 2cell, dvips]{xy}
\usepackage{latexsym,amsmath, amsfonts, graphics, epsf, epic}
\usepackage{amsthm}
\usepackage{amssymb}
\usepackage{amsopn}
\usepackage{amscd}
\usepackage{color}

\newtheorem{thm}{Theorem}[section]
\newtheorem{cor}[thm]{Corollary}
\newtheorem{lem}[thm]{Lemma}

\newtheorem{defn}[thm]{Definition}

\DeclareMathOperator{\Hom}{Hom} \DeclareMathOperator{\End}{End}
\DeclareMathOperator{\sgn}{sgn} 
\DeclareMathOperator{\M}{\textsf{M}}
\DeclareMathOperator{\PW}{P^+(n,r)}
\DeclareMathOperator{\LW}{\Lambda^+}
\DeclareMathOperator{\W}{P^+(n,r)}

\newcommand{\To}{\longrightarrow}

\DeclareMathOperator{\he}{\hat{e}}
\DeclareMathOperator{\hv}{\hat{v}}

\newcommand{\dn}{\underline{n}}
\newcommand{\di}{\underline{i}}
\newcommand{\kk}{\underline{k}}
\newcommand{\db}{\underline{b}}
\newcommand{\da}{\underline{a}}
\newcommand{\jj}{\underline{j}}
\newcommand{\dk}{\underline{k}}
\newcommand{\dl}{\underline{l}}

\newenvironment{pf}{\noindent \textit{Proof.}\,}{\hfill$\square$ \vskip5pt}

\title{Comparing $GL_n$-Representations by Characteristic-Free
Isomorphisms between Generalized Schur Algebras}

\author{Ming Fang, Anne Henke and Steffen Koenig} 

\thanks{Ming Fang acknowledges support from AsiaLink grant
ASI/B7-301/98/679-11.}

\date{\today}

\begin{document}

\maketitle

\centerline{With an appendix by {\sc Stephen Donkin}}

\begin{abstract}
Isomorphisms are constructed between generalized Schur algebras in
different degrees. The construction covers both the classical case
(of general linear groups over infinite fields of arbitrary
characteristic) and the quantized case (in type $A$, for any
non-zero value of the quantum parameter $q$). The construction
does not depend on the characteristic of the underlying field or
the choice of $q \neq 0$. The proof combines a combinatorial
construction with comodule structures and Ringel
duality. Applications range from equivalences of categories to results
on the structure and cohomology of Schur algebras to identities of
decomposition numbers and also of $p$-Kostka numbers, in both cases
reproving and generalizing row and column removal rules.
\end{abstract}

\markright{\sc Comparing Schur Algebras}

\section{Introduction}
\subsection{General linear groups and Schur algebras}
Let $k$ be an infinite field of arbitrary characteristic.  The
homogeneous polynomial representations of the general linear group
$GL_n(k)$ of a fixed degree $r$ are precisely the modules over
(classical) Schur algebras $S_k(n,r)$.  Classical and quantized Schur
algebras and their representations have been used in a variety of
contexts, both on structural and on numerical level. Among these is
the representation theory of the algebraic group $GL_n$ and of the
finite group $GL_n$ in describing and in cross characteristic,
representation theory of symmetric groups and of their Hecke algebras,
polynomial functors and group cohomology.
This article deals with classical Schur algebras, whose representations
are precisely the polynomial representations of the algebraic group
$GL_n(k)$ over any infinite field $k$; simultaneously we deal with 
the quantized Schur algebras.

\subsection{The main result}
Let $GL_n=GL_n(k)$ for a fixed infinite field $k$.  
The simple $GL_n$-representations
of degree $r$ are parametrized by partitions $\lambda \in
\Lambda^+(n,r)$ of $r$ into not more than $n$ parts.
A {\it saturated} subset of $\Lambda^+(n,r)$ is a subset $\Pi$ such
that $x \in \Pi$ whenever $x < y \in \Pi$.  Choose any saturated
subset $\Pi \subset \Lambda^+(n,r)$. Moreover, choose any $m \in
\mathbb{N}$ such that the rows of the partitions in $\Pi$ have not
more than $m$ boxes, that is, each partition $\lambda \in \Pi$ fits
into an $n \times m$ rectangle. Let $\hat{\lambda}$ be the complement
of $\lambda$ within this rectangle; then $\hat{\lambda}$ is a
partition of $mn-r$.  Our Main Theorem \ref{maintheorem} implies
equivalences of categories as follows:

{\em The abelian category of
  $GL_n$-modules of degree $r$ with simple composition factors indexed
  by $\lambda \in \Pi$ is equivalent to the category of $GL_n$-modules
  of degree $mn-r$ with simple composition factors
indexed by $\hat{\lambda}$.} 
This equivalence implies identities between $GL_n$-decomposition
numbers in degrees $r$ and $nm-r$, and isomorphisms in cohomology
between representations in degrees $r$ and $nm-r$. The equivalence
preserves vector space dimensions.

Stronger and more explicitly, 
the {\em Main Theorem} \ref{maintheorem} states the
existence of an isomorphism between quotient algebras of Schur
algebras (modulo ideals in heredity chains) as follows: {\em For any
  natural numbers $n$ and $r$ and for any saturated set $\Pi\subset
  \Lambda^+(n,r)$, and for any natural number $m$ such that the
  partitions in $\Pi$ have rows of length not more than $m$, there is
  an isomorphism between (classical or quantized) generalized Schur
  algebras:
  $$S(n,r)/I_{\Pi} \simeq S(n,nm-r)/I_{\widehat{\Pi}}.$$} 
Here,
taking complements of partitions defines a bijection between the index
sets $\Pi$ and $\widehat{\Pi}$, which determine the ideals $I_{\Pi}$
and $I_{\widehat{\Pi}}$ in the heredity chain. The isomorphisms of
algebras will be {\em constructed explicitly}.

Quotients as in the main theorem sometimes are called generalized
Schur algebras.  Behind these isomorphisms is another, more general,
set of isomorphisms on centralizer subalgebras of Ringel duals of
Schur algebras.

The above statement about equivalences of categories follows from
these isomorphisms. The correspondence between projective modules of
the two algebras is again given by sending a partition $\lambda$ to
its complement $\hat{\lambda}$, defined above.

The Main Theorem \ref{maintheorem} covers not only the case of
classical Schur algebras, but also that of quantized Schur algebras
(again for general linear groups). In this case, Theorem
\ref{maintheorem} generalizes a result of Beilinson, Lusztig and
MacPherson \cite{BLM}, who constructed isomorphisms between a
quantized Schur algebra $S(n,r)$ and a quotient of $S(n,r+n)$, in
order to describe quantized enveloping algebras geometrically (by
transfering a description via convolution of 'pairs of flags' from
quantized Schur algebras to an inverse limit of $q$-Schur algebras,
which contains the quantized enveloping algebra). These isomorphisms
can be recovered as products of two maps in the family of
isomorphisms provided by Theorem \ref{maintheorem}.

\subsection{Outline of the article}
For background on Schur algebras, the reader is refered to
Green \cite{Green} (classical case) and Donkin \cite{Donkin}
(quantized case). Additional information and references can be found
in \cite{K}.  Here is an outline of the article: In Section 2 we
collect basic information about Schur algebras, and we fix some
notation. Sections 3 to 6 are devoted to the proof of the Main Theorem
\ref{maintheorem}. The ingredients of the proof are:
\begin{itemize}
\item[(a)] combinatorics of a construction sending
a partition $\lambda$ to its 'complement' $\hat{\lambda}$ (Section 3);
\item[(b)] a result (Theorem \ref{compareGmaps}), identifying
homorphisms of certain representations (tensor products of exterior
powers) of classical or quantized $GL_n$ in different degrees by using
equations given by comodule structures (Section 4);
\item[(c)] Ringel duality for quasi-hereditary algebras (Section
5), using that the representations occuring in Section 5 are so-called
tilting modules.
\end{itemize}
In Section 6, the Main Theorem is stated in full generality and its
proof is completed by putting together the earlier results. In Sections
7 to 9 we collect some applications to decomposition numbers, to cohomology,
to group algebras of symmetric groups and Hecke algebras and
in particular to $p$-Kostka numbers. We also
recover the isomorphisms constructed in \cite{BLM}, implying
James' column removal rule. Moreover, we construct
a Morita equivalence between centralizer subalgebras of generalized
Schur algebras, which implies James' row removal rule.

The results of this article, which are characteristic independent,
complement the characteristic dependent isomorphisms of centralizer
subalgebras of (classical) Schur algebras constructed in \cite{H} and
also in \cite{HK}.

The appendix, written by Stephen Donkin, provides an alternative 
approach to and a generalization of the Main Theorem.

\section{Setup and notation}

\subsection{Notation} \label{21}

Throughout, let $k$ be an infinite field.  For a non-zero natural
number $n$ define $\dn = \{1,\ldots,n\}$. Let $r$ be another non-zero
natural number. We define sets of multi-indices of length $r$ with
entries in $\underline{n}$ as follows:
\begin{eqnarray*}
I(n,r) &=& \{ \di =(i_1,\ldots,i_r) \mid i_\rho\in \dn \}, \\
I^{-}(n,r) &=& \{ \di \in I(n,r) \mid  i_1> \ldots > i_r \}.
\end{eqnarray*}
Denote by $\Lambda^+(n,r)$ the set of partitions of $r$ with at
most $n$ parts and by $P(n,r)$ the set of sequences of
non-negative integers $(a_1,a_2, \ldots)$ such that $0\leq
a_i\leq n$ and $a_1+a_2+\ldots = r$.  Let $\PW$ be the subset of
partitions in $P(n,r)$, and $\PW_m$ the subset of $\PW$ such that
$a_{m+1}= a_{m+2} =\ldots = 0$.
The ordering considered on both $\Lambda^+(n,r)$ and $\PW$ is the dominance
ordering, namely,
\[
\lambda \geq \mu\iff \lambda_1+\ldots
+\lambda_s\geq\mu_1+\ldots+\mu_s
\]
for all natural numbers $s \geq 1$ and $\mid \lambda \mid \,=\,\mid
\mu \mid$.  Given a partition $\alpha\in P^+(n,r)$, we define
the following sets of tableaux:
\begin{eqnarray*}
Tab(\alpha)&=&\{ \alpha'-\mbox{tableaux with entries in $\dn$} \}, \\
Tab^+(\alpha)&=& \{ S \in Tab(\alpha) \mid
\mbox{ entries are strictly increasing down columns} \}, \\
Tab^-(\alpha)&=& \{ S \in Tab(\alpha) \mid
\mbox{  entries are strictly decreasing down columns} \}, \\
STab(\alpha)&=& \{ S \in Tab^-(\alpha) \mid \mbox{ entries are
weakly decreasing along its rows} \}.
\end{eqnarray*}
The set $STab$ is in one-to-one correspondence with the set of all
semi-standard tableaux $S \in Tab(\alpha)$, that is by definition the
set of all tableaux $S$ whose entries are strictly increasing down
columns and weakly increasing along its rows.  By abuse of notation,
we call a tableau in $STab$ also semi-standard.  Note that we are
always using the conjugate partition $\alpha'$ in these definitions.

\subsection{Schur algebras}
\label{sec1}
Let $G=GL_n$ and denote by $k^{G}$ the set of all the $k$-valued
functions on $G$. Define maps $$\Delta : k^G\to k^{G\times G}
\mbox{ and } \varepsilon : k^G\to k$$ by setting $\Delta(f)(g_1,
g_2)=f(g_1g_2)$ and $\varepsilon (f) = f(1)$ for any $f\in k^G$ and
$g_1, g_2\in G$.  Here $1$ denotes the identity matrix of size
$n\times n$.  We will identify
$k^G \otimes k^G$ with a subspace of $k^{G \times
G}$ by $$(f_1 \otimes f_2)(g,h)= f_1(g)f_2(h).$$

For each pair $(\mu,\nu) \in \dn \times \dn$,
let $c_{\mu\nu}$ be the function
which associates to each $g\in G$ its $(\mu,\nu)$-entry
$g_{\mu\nu}$. Denote by $A_k(n)$ the $k$-subalgebra of $k^G$
generated by $c_{\mu\nu}\ (\mu, \nu\in \dn)$. Since $k$ is
infinite, this is precisely the polynomial ring over $k$ in $n^2$
indeterminates $c_{\mu\nu}$. Let $A_k(n, r)$ be the homogenous
subspace in $A_k(n)$ of degree $r$ in the indeterminates
$c_{\mu\nu}$.  Then $A_k(n, r)$ as a $k$-vector space is spanned
by
\[
c_{\di,\jj} = c_{i_1j_1}c_{i_2j_2}\cdots c_{i_rj_r},
\]
for all $\di,\jj \in I(n,r)$. As in Green~\cite{Green}, we separate in
this notation multi-indices $\di, \jj$ by a comma, while we do not
separate natural numbers $i_1 j_1$.  The maps
$\Delta$ and $\varepsilon$ defined above, when applied to the functions
$c_{\di,\jj}$, satisfy the following equations:
\begin{equation*}
  \Delta (c_{\di,\jj}) = \sum_{\kk\in I(n,r)} c_{\di,\kk}\otimes c_{\kk,\jj},
  \qquad \varepsilon (c_{\di,\jj}) = \delta_{\di,\jj}
\end{equation*}
where $\delta_{\di,\jj} = 1$ or $0$, if $\di=\jj$ or $\di \neq \jj$.
The vector space $A_k(n,r)$ forms a coalgebra with comultiplication
$\Delta$ and counit $\varepsilon$.  The (classical) Schur algebra,
denoted by $S_k(n,r)$, is the dual of the coalgebra $A_k(n,r)$.  In
the following we will omit the subscript $k$ (to indicate the
underlying field) on coalgebras and algebras.

Let $W$ be a representation of $G$.  Suppose $W$ has basis
$\{w_a\}$.  Then we call $W$ a {\it homogeneous polynomial
representation of degree $r$} (in defining characteristic) if
\[
g\cdot w_b = \sum _a r_{ab}(g)w_a
\]
where all $r_{ab}\in A(n,r)$.  By \cite{Green}, homogeneous
polynomial representations
of degree $r$ are precisely the modules for the
Schur algebra $S(n,r)$.  Alternatively, they are precisely the
comodules for the coalgebra $A(n,r)$. Note that the natural
representation of $G$ is polynomial of degree one. Moreover,
tensor products, symmetric powers and exterior powers of
polynomial representations are again polynomial. More details
about exterior powers are given in Section~\ref{sec2} where their
structure maps as comodules are given.

In the quantized case, the Schur algebra is defined in a similar way.
For details, the reader is refered to Donkin~\cite{Donkin}.  In this
paper, groups and algebras are always considered over a field $k$, and
the quantum parameter $q \in k$ always will be assumed to be non-zero;
in particular this implies that $q$ is invertible.

\subsection{Exterior powers} \label{sec2}

For non-zero natural numbers $n$ and $r$, denote by $S(n,r)$
either a classical Schur algebra as discussed by Green in
\cite{Green} for the algebraic group $GL_n$ (see Section
\ref{sec1}), or a $q$-Schur algebra as defined in Donkin's book
\cite{Donkin}. We assume that the quantum parameter $q$ is
different from zero. Let $\M$ be the classical or quantum monoid
which affords precisely the polynomial representations of
classical or quantum $GL_n$ respectively.
By $E$ or $V$ we denote the natural left or right $\M$-module
respectively -- that means right or left $k[\M]$-comodule; note
the switch of sides of the actions when translating from
$\M$-modules to $k[\M]$-comodules. Let $E$ have $k$-basis
$e_1,\ldots, e_n$ and suppose $V$ has $k$-basis $v_1,\ldots,
v_n$.  Viewed as a comodule over the coordinate ring $k[\M]$,
then $E$ and $V$ have structure maps:
\begin{eqnarray*}
 \tau_{1,E} (e_i)= \sum_{j=1}^n e_j\otimes c_{ji},  & &
\tau_{1,V}(v_i) = \sum_{j=1}^n c_{ij}\otimes v_j
\end{eqnarray*}
where $1 \leq i \leq n$.

The exterior powers of $E$ and $V$ are again left and right
$\M$-modules. Using the notation for sets of multi-indices introduced
in Section \ref{21}, as $q \neq0$, the exterior power
$\wedge^rV$ has $k$-basis $\{\hat{v}_{\di} \mid \di\in I^-(n,r)\}$ and
$\wedge^rE$ has $k$-basis $\{\hat{e}_{\di} \mid \di\in I^-(n,r)\}$
where $\hat{v}_{\di}= v_{i_1} \wedge \ldots \wedge v_{i_r}$ and
$\hat{e}_{\di}= e_{i_1} \wedge \ldots \wedge e_{i_r}$.  Note that
$\wedge^rV=0=\wedge^rE$ if $r > n$.  The structure maps of the
exterior powers are given by (see \cite{Donkin}, 1.3.1):
\begin{eqnarray*}
\tau_{r,E} (\hat{e}_{\di})= \sum_{\jj \in I^-(n,r)} \hat{e}_{\jj}\otimes
\langle \jj:\di \rangle, & & \tau_{r,V} (\hat{v}_{\di}) = \sum_{\jj \in
I^-(n,r)} (\di:\jj) \otimes \hat{v}_{\jj},
\end{eqnarray*}
with bideterminants \label{definestructuremaps}
\begin{eqnarray} \label{eqn1}
\langle \di:\jj \rangle =
\sum_{\pi \in \Sigma_r} (-q)^{l(\pi)} c_{\di\pi,\jj}, & &
(\di:\jj) = \sum_{\pi \in \Sigma_r} (-1)^{l(\pi)} c_{\di,\jj \pi}.
\end{eqnarray}
Here $l(\pi)$ denotes the length of the permutation $\pi$.
For any $\alpha\in \PW$ and $S\in Tab(\alpha)$, define
\begin{eqnarray*}
\he_S &=& e_{S(1,1)}\wedge\ldots\wedge e_{S(\alpha_1,1)}
\otimes e_{S(1,2)}\wedge\ldots\wedge e_{S(\alpha_2,2)}\otimes \ldots, \\
\hv_S &=& v_{S(1,1)}\wedge\ldots\wedge v_{S(\alpha_1,1)}
\otimes v_{S(1,2)}\wedge\ldots\wedge v_{S(\alpha_2,2)}\otimes \ldots.
\end{eqnarray*}
Here $S(i,j)$ denotes the entry in row $i$ and column $j$ of tableau
$S$.  Note that $\hv_S=0=\he_S$ if $S$ has a repetition in the entries
of one of its columns. The tensor product of the exterior powers
$\wedge^{\alpha_i}V$, denoted by $\wedge^{\alpha}V$, has $k$-basis
$\{\hv_S : S\in Tab^-(\alpha)\}$, and $\wedge^{\alpha}E$ has $k$-basis
$\{\he_S : S\in Tab^-(\alpha)\}$.  We denote by $\tau_{\alpha,E}$ and
$\tau_{\alpha,V}$ the structure maps of $\wedge^{\alpha}E$ and
$\wedge^{\alpha}V$ respectively; we write shortly $\tau_\alpha$ when
no confusion arises about the vector space used. Define
bideterminants
\begin{eqnarray*}
(S:T)=\prod_{}^{} (S_i:T_i) & \mbox{ and }&
\langle S:T \rangle=\prod_{}^{} \langle S_i:T_i \rangle
\end{eqnarray*}
Here we consider the columns (read from top to bottom) of a
tableau $S$ as multi-indices and write $S_i$ for the $i$-th
column of the tableau $S$. Note that, compared to Lemma 1.3.1 in
\cite{Donkin}, there is a switch of rows and columns here; this
is due to the different definition of the set $Tab^-$.

For later use, we next
collect a few technical results; in particular, property (3) of the
following Lemma can be taken as the definition of the quantum
determinant.
\begin{lem} \label{lem2} \label{lem1}
With the notation as above, the following statements hold true:
\begin{enumerate}
\item Let $S \in Tab^-(\alpha)$, then the structure maps
  $\tau_{\alpha,E}$ and $\tau_{\alpha,V}$ are given by
\begin{eqnarray*}
& & \tau_{\alpha,E} (\he_S) = \sum_{T\in Tab^-(\alpha)}
\he_T\otimes \langle T:S \rangle,  \\ & & \tau_{\alpha,V} (\hv_S)
= \sum_{T\in Tab^-(\alpha)} (S:T)\otimes \hv_T.
\end{eqnarray*}
\item For any $S, T\in Tab^-(\alpha)$ we have
$(S:T) = \langle S:T \rangle$.
\item We have $\det= (\di:\di) = \langle \dl : \dl \rangle$, where
  $\di=(1,2,\ldots,n)$ and $\dl=(n,\ldots, 2,1)$. In particular, if
  $\alpha=(n, \ldots , n)\in \Lambda^+(n,nm)$ then $\tau_{\alpha}$ is
  the map given by multiplication with $(\det)^m$.
\end{enumerate}
\end{lem}
\begin{pf}
  The first two statements can be found in \cite{Donkin}, Lemma 1.3.1,
  the third statement in \cite{DD}, Theorem 4.1.7 (see also
  \cite{Donkin}, Section 0.20). In case
  $\alpha=(n, \ldots , n)$, note that there is only one element, say
  $A$, in $Tab^-(\alpha)$, namely the tableau whose columns are all
  filled decreasingly with $n$ to $1$. In this case $(A:A)=
  \prod_{i=1}^m (\dl:\dl)$ with $\dl$ as above.
\end{pf}

\section{Combinatorics}\label{sec3}

\subsection{Complements of Partitions.} Given a partition
$\lambda = (\lambda_1,\ldots,\lambda_n)$ in $\Lambda^+(n,r)$ and a
positive integer $m \geq \lambda_1$, we define the {\it complement
  partition} $T_{m}(\lambda)$ or $\hat{\lambda}$ of $\lambda$ with respect
to $m$ as follows:
\[
T_{m}(\lambda) = \hat{\lambda} = (m - \lambda_n, m -
\lambda_{n-1}\ldots,m-\lambda_{1}).
\]
So the function $T_{m}$ takes the complement of $\lambda$ within
an $n \times m$-rectangle. Note that it really depends from both
$m$ and $n$. Since we only apply the function when $n$ is fixed,
we write shortly $T_m$. When both $n$ and $m$ are understood, we
write $\hat{\lambda}$.
For example, let $\lambda= (4,2,1)$ in $\Lambda^+(4,7)$, so $n=4$.
Then $T_4(\lambda) = (4 - \lambda_4, 4 - \lambda_3, 4 - \lambda_2,
4 - \lambda_1) = (4,3,2,0)$ whereas
$T_5(\lambda) = (5,4,3,1).$ Observe that for $m=5 \geq \lambda_1=4$
the coloured partition in the picture is $T_5(\lambda)$; the remaining
white boxes of the $n \times m$ rectangle are a skew-diagram
$\iota(\lambda)$, obtained by reflecting $\lambda$ at its top left corner.

$$
\vspace{1.3cm}
\setlength{\unitlength}{1.1mm}
\begin{picture}(-120,0)(50,11)

\thinlines
\drawline[20](0,0)(0,16)
\drawline[20](0,16)(20,16)
\drawline[20](20,16)(20,0)
\drawline[20](20,0)(0,0)
\drawline[20](4,0)(4,16)
\drawline[20](8,0)(8,16)
\drawline[20](12,0)(12,16)
\drawline[20](16,0)(16,16)
\drawline[20](0,4)(20,4)
\drawline[20](0,8)(20,8)
\drawline[20](0,12)(20,12)
\put(-35,12){\makebox(0,0){$(5,4,3,1)=T_5(4,2,1)=T_{m}(\lambda) 
\longrightarrow$}}\put(35,4){\makebox(0,0){$\longleftarrow \iota(\lambda)$}}

\textcolor{yellow}
{\linethickness{4mm}
\put(2,14){\picsquare}
\put(6,14){\picsquare}
\put(10,14){\picsquare}
\put(14,14){\picsquare}
\put(18,14){\picsquare}
\put(2,10){\picsquare}
\put(6,10){\picsquare}
\put(10,10){\picsquare}
\put(14,10){\picsquare}
\put(2,6){\picsquare}
\put(6,6){\picsquare}
\put(10,6){\picsquare}
\put(2,2){\picsquare}}
\end{picture}
$$

In later sections of this paper, we will also work with the
conjugate situation: instead of taking a complement of a
partition $\lambda \in \Lambda^+(n,r)$ with $\lambda_1 \leq m$ in
an $n \times m$-rectangle, we take the complement of the
conjugate partition $\lambda'=\alpha \in P^+(n,r)_m$ in an $m
\times n$-rectangle; in the situations occuring, $m$ is fixed and
we write $n-\alpha$ for the complement in the conjugate
situation; in particular $\hat{\lambda} = T_m(\lambda)
=(n-\alpha)'$.

We collect some simple properties of forming the complement of a partition.
The proof is left to the reader.

\begin{lem} \label{prop1}
  For partitions $\lambda,\mu\in \Lambda^+(n,r)$ and a natural number $m\geq
  \lambda_1$ we have:
$T_{m}^2(\lambda) =\lambda$, and $\lambda\geq \mu$ implies
$T_{m}(\lambda)\geq T_{m}(\mu)$.
\end{lem}

\subsection{Complements of tableaux.}
If $n \geq r$, then for any multi-index $\di$ in $I^-(n,r)$, we
define the {\it complement} $\hat{\di}$ of $\di$ by $\hat{\di} =
(i_{r+1},\ldots, i_n)$ with $i_{r+1}> \ldots> i_n$  such that
$\di\cup \hat{\di} = \dn$ and $\di\cap \hat{\di}=\emptyset$;
here, by abuse of notation, we consider multi-indices $\jj=(j_1,
\ldots, j_r)$ of some length $r$ with pairwise different entries
as sets with elements $j_\rho$, for $1 \leq \rho \leq r$. As
before, we consider the columns, read from top to bottom, of a
tableau $S$ as multi-indices. If $S_i$ is the $i$-th column of
the tableau $S$ and has no repeated entries, then $\widehat{S_i}$
is the complement of the multi-index $S_i$, as defined above.

For a partition $\lambda$ with $m\geq
\lambda_1$, and a $\lambda$-tableau $S$ with entries from the set
$\underline{n}$, whose entries are strictly decreasing down
columns, we shall define the complement $T_{m}(\lambda)$-tableau
$T_{m}(S)$ by
\[
T_{m}(S)_i = \widehat{S_{m+1-i}}.
\]
If $m$ is fixed, we will write $\widehat{S}$ for $T_{m}(S)$. For
example, let $\lambda= (4,2,1)$ in $\Lambda^+(4,7)$ and $m=5$ as
before and let $S$, filled with entries in $\{1, \ldots, 4 \}$,
be given by

$$
\vspace{.3cm}
\setlength{\unitlength}{1.1mm}
\begin{picture}(-50,0)(50,3)
\thinlines
\drawline[20](0,12)(16,12)
\drawline[20](0,8)(16,8)
\drawline[20](0,4)(8,4)
\drawline[20](0,0)(4,0)
\drawline[20](0,0)(0,12)
\drawline[20](4,0)(4,12)
\drawline[20](8,4)(8,12)
\drawline[20](12,8)(12,12)
\drawline[20](16,8)(16,12)
\put(19,0){.}
\put(1.25,8.75){4} \put(5.25,8.75){3} \put(9.25,8.75){3}
\put(13.25,8.75){2}
\put(1.25,4.75){3} \put(5.25,4.75){2}
\put(1.25,0.75){1}
\end{picture}
$$

Define the skew-tableau $\iota(S)$ by reflecting the tableau $S$ at
its top left corner.  Then $T_5(S)$ is obtained by placing $\iota(S)$
in the right-hand bottom corner of an $n=4$ times $m=5$ box and by
then filling the complement of each column into the empty boxes of
that column, for example:

$$
\vspace{1.3cm}
\setlength{\unitlength}{1.1mm}
\begin{picture}(-50,0)(50,11)
\thinlines
\drawline[20](0,0)(0,16)
\drawline[20](0,16)(20,16)
\drawline[20](20,16)(20,0)
\drawline[20](20,0)(0,0)
\drawline[20](4,0)(4,16)
\drawline[20](8,0)(8,16)
\drawline[20](12,0)(12,16)
\drawline[20](16,0)(16,16)
\drawline[20](0,4)(20,4)
\drawline[20](0,8)(20,8)
\drawline[20](0,12)(20,12)
\put(-10,9){\makebox(0,0){$T_5(S)=$}}
\put(35,4){\makebox(0,0){with $\iota(S)=$}}
\put(1.25,12.75){4} \put(5.25,12.75){4} \put(9.25,12.75){4}
\put(13.25,12.75){4} \put(17.25,12.75){2}
\put(1.25,8.75){3} \put(5.25,8.75){3} \put(9.25,8.75){2}
\put(13.25,8.75){1}
\put(1.25,4.75){2} \put(5.25,4.75){1} \put(9.25,4.75){1}
\put(1.25,0.75){1}
\drawline[20](50,0)(66,0)
\drawline[20](50,4)(66,4)
\drawline[20](58,8)(66,8)
\drawline[20](62,12)(66,12)
\drawline[20](50,0)(50,4)
\drawline[20](54,0)(54,4)
\drawline[20](58,0)(58,8)
\drawline[20](62,0)(62,12)
\drawline[20](66,0)(66,12)
\put(51.25,0.75){2} \put(55.25,0.75){3} \put(59.25,0.75){3}
\put(59.25,4.75){2} \put(63.25,0.75){4} \put(63.25,4.75){3}
\put(63.25,8.75){1} \put(68,0){.}
\end{picture}
$$

\begin{thm}\label{ssmirror}
  For any $\lambda\in \Lambda^+(n,r)$ and integer $m\geq \lambda_1$,
  $S$ is a semi-standard $\lambda$-tableau if and only if $T_{m}(S)$ is
  a semi-standard $T_{m}(\lambda)$-tableau.
\end{thm}
\begin{pf}
  Without loss of generality, it is enough to prove the claim in the
  case where $\lambda$ has only two columns and ${ m}=2$.
  We need a partial ordering defined on the set of all multi-indices.
  Namely, for $\di=(i_1,\ldots,i_v)\in I(n,v)$ and $\jj=
  (j_1,\ldots,j_w)\in I(n,w)$, write $\jj \preceq \di$ if $w\leq v$
  and $j_1\leq i_1,\ldots, j_w \leq i_w$.

  Let $\lambda$ be a two-column partition in $\Lambda^+(n,r)$.
  Consider a semi-standard tableau $S$ of shape $\lambda$ with columns
  $S_1$ and $S_2$, filled by entries in $\underline{n}$.
Since $S$ is semi-standard, both its columns can be considered as
multi-indices, say
\[
S_1 = \di = (i_1,\ldots,i_v) \mbox{\hspace{.2cm} and \hspace{.2cm}}
S_2 = \jj = (j_1,\ldots,j_w)
\]
where $\di$ and $\jj$ are ordered by size and each have entries in
$\underline{n}$ without repetitions. By the definition in Section 2.1,
$S \in STab$ means that entries in rows of $S$ are weakly decreasing
and entries in columns of $S$ are strictly decreasing.  Since $S$ is
semi-standard, then $w \leq v$ and by definition $\jj\preceq\di$.

Let $\da$ and $\db$ be any multi-indices and define the operation
$P_{i,j}$ on the multi-indices $\da$ and $\db$ for $i,j\in
\{0,1,\ldots,n\}$ as follows: $P_{i,j}$ bumps out $i$ from $\da$ and
$j$ from $\db$.  If $i\notin \da$ (resp.~ $j\notin \db$), then
$P_{i,j}$ does nothing to $\da$ (resp.~ $\db$).
For example, if $i = a_s$ and $j = b_t$, then the $\rho$th entries of
the multi-indices $P_{i,j}\da$ and $P_{i,j}\db$ are
\[
[P_{i,j}\da]_\rho =\left\{
\begin{array}{ll}
a_\rho & \rho<s,\\
a_{\rho+1} & \rho\geq s,
\end{array}
\right.\qquad
[P_{i,j}\db]_\rho =\left\{
\begin{array}{ll}
b_\rho & \rho<t,\\
b_{\rho+1} & \rho\geq t.
\end{array}
\right.
\]
The length of the multi-indices $P_{i,j}\da$ and $P_{i,j}\db$ are now
both one shorter than that of $\da$ and $\db$ respectively.
$$
\vspace{1.7cm}
\setlength{\unitlength}{1.4mm}
\begin{picture}(-90,0)(50,13)

\thinlines
\drawline[20](0,2)(0,14)
\drawline[20](4,2)(4,14)
\drawline[20](8,8)(8,14)
\drawline[20](0,14)(8,14)
\drawline[20](4,8)(8,8)
\drawline[20](0,2)(4,2)
\put(-6,8){\makebox(0,0){$S=$}}
\put(2,11){\makebox(0,0){$\di$}}
\put(6,11){\makebox(0,0){$\jj$}}
\drawline[20](37,2)(37,14)
\drawline[20](41,2)(41,14)
\drawline[20](45,8)(45,14)
\drawline[20](37,14)(45,14)
\drawline[20](41,8)(45,8)
\drawline[20](37,2)(41,2)
\put(24,8){\makebox(0,0){and $T_2(S)=\widehat{S}=$}}
\put(39,11){\makebox(0,0){$\hat{\jj}$}}
\put(43,11){\makebox(0,0){$\hat{\di}$}}
\end{picture}
$$
We now fix $\da = (n, \dots, 1) = \db$ as multi-indices containing the
numbers $1, \dots, n$ in strictly decreasing order.  With the
operations $P_{i,j}$ we can write $\hat{\di}$ and $\hat{\jj}$ as
follows:
\begin{eqnarray} \label{eqn77}
\hat{\di} &=& P_{i_v,0}\cdots P_{i_{w+1},0}\cdot P_{i_w,j_w}
\cdots P_{i_1,j_1}(\da), \\
\label{eqn78} \hat{\jj} &=& P_{i_w,j_w}\cdots P_{i_1,j_1}(\db).
\end{eqnarray}
Note that the two columns of $\widehat{S}=T_2(S)$ are exactly $\hat{\jj}$
and $\hat{\di}$. So to prove $T_2(S)$ is semi-standard, we only need
to show that $\hat{\di}\preceq\hat{\jj}$.
We do this by induction on the number of the above operations
$P_{i,j}$ that we take. More precisely, we start with the
multi-indices $\da = (n, \dots, 1) = \db$, which certainly satisfy
$\da \preceq \db$. Then we successively apply operations $P_{i,j}$
with pairs of indices $(i,j)$ where $i$ and $j$ are from the sets
${\di}$ and ${\jj} \cup \{ 0 \}$ respectively. In all $P_{i,j}$
used in the above Equations (\ref{eqn77}) and (\ref{eqn78}), we
always have $i \geq j$. And we never repeat any index $i$ or $j$
except possibly $j=0$. Suppose after $\rho\leq w$ such
operations, we get multi-indices $\dk$ and $\dl$ with $\dk
\preceq \dl$.  Then for $p\in \di$, $q\in \jj$ with $p\geq q$, we
must have $p = k_s, q = l_t$ with $s \leq t$. Therefore,
\begin{align*}
[P_{p,q}(\dk)]_\rho = [\dk]_\rho, & \qquad 1\leq \rho\leq s-1,\\
[P_{p,q}(\dl)]_\rho = [\dl]_\rho, & \qquad  1\leq \rho\leq t-1.
\end{align*}
Hence,
\begin{align*}
[P_{p,q}(\dk)]_\rho &\leq [P_{pq}(\dl)]_\rho, & 1\leq \rho \leq s-1,\\
[P_{p,q}(\dk)]_\rho &=[\dk]_{\rho+1}\leq [\dl]_{\rho+1}\leq [\dl]_\rho =
[P_{p,q}(\dl)]_\rho, & s\leq \rho\leq t-1,\\
[P_{p,q}(\dk)]_\rho &=[\dk]_{\rho+1}\leq [\dl]_{\rho+1} =
[P_{p,q}(\dl)]_\rho, & t\leq \rho\leq n.\\
\end{align*}
Hence $P_{p,q}(\dk)\preceq P_{p,q}(\dl)$, which is also true in case
$q=0$. Inductively,
\[
P_{i_w,j_w}\cdots P_{i_1,j_1}(\da)\preceq P_{i_w,j_w}\cdots
P_{i_1,j_1}(\db)
\]
which implies immediately that
\[
\hat{\di} =
P_{i_v,0}\cdots P_{i_{w+1},0}\cdot P_{i_w,j_w}\cdots P_{i_1,j_1}(\da)
\preceq P_{i_w,j_w}\cdots P_{i_1,j_1}(\db)   =\hat{\jj}
\]
Since the property of being semistandard is determined by comparing any
two neighbouring columns, the proof generalises to
any semi-standard tableau $S$ of shape $\lambda \in \Lambda^+(n,r)$
which is filled by entries in $\underline{n}$.
\end{pf}

\subsection{Signs of tableaux}
Using earlier notations, we finally introduce the sign of a
tableau, needed for the explicit construction of the isomorphisms
between quantized Schur algebras in the next section. In Section
3.2, we defined for $n\geq r$ and any multi-index $\di\in
I^-(n,r)$ its complement $\hat{\di}\in I^-(n, n-r)$ to be the
complement of $\di$ in $\dn$.
Let $\omega = v_n\wedge\ldots\wedge v_1$, a
generator of $\wedge^nV$, and define the sign of $\di\in I^-(n,r)$
by the equation
\[
v_{\di}\wedge v_{\hat{\di}}  = \sgn(\di)\omega.
\]
We generalize the definition of a sign for a multi-index to that
of a sign of a tableau $S\in Tab^-(\alpha)$ with $m$ columns $S_1,
\dots, S_m$ as follows:
\[
\sgn(S) = \sgn(S_1) \cdot \sgn(S_2) \cdots \sgn(S_m)
\]
Moreover, we define
\[
|\di| = i_1 + \ldots + i_r
\]
and
\[
\epsilon(S) = (m-1)|\widehat{S_1}|+(m-2)|\widehat{S_2}|+\ldots + 2
|\widehat{S_{m-2}}| + |\widehat{S_{m-1}}|,
\]
where $\widehat{S_i}$ is the complement of the multi-index $S_i$.

\section{Homomorphisms between exterior powers}
Fix a natural number $m$.  Then for any $\alpha =
(\alpha_1,\ldots,\alpha_m)\in \W_m$, we denote by $n-\alpha$ the
partition $(n- \alpha_m,\ldots,n-\alpha_1)$.  Note that this depends
on the choice of $m$, compare with $T_{m}$ in Section 3.
In this section, we compare morphism spaces between tensor products of
exterior powers $\wedge^i$ and $\wedge^{n-i}$ respectively. The main
result identifies such two morphism spaces in different degrees.  We
will use the notation for exterior powers introduced in Section
\ref{sec2}. In particular, $V$ is an $n$-dimensional vector space.
Then we have for $1 \leq r \leq n$ the following result -- which is a
special case of~\cite{Donkin}, Lemma 1.3.3:

\begin{lem} \label{lem3}
  The linear map $\mu=\mu_r: \wedge^nV\otimes \wedge^rV\To
  \wedge^rV\otimes \wedge^nV$ defined by
\[
\mu(\omega\otimes \hv_{\di}) = q^{-|\di|}\hv_{\di}\otimes \omega
\]
for $\di\in I^-(n,r)$ is an $\M$-module homomorphism.
\end{lem}

\begin{pf}
  By \cite{DD} Theorem 4.1.9 or \cite{Donkin} Section 0.20, we have
  $q^j \cdot c_{ij} \cdot \det = q^{i} \cdot \det \cdot c_{ij}$ for $1
  \leq i,j \leq n$, where $det = (\dl:\dl)$ with $\dl=(n, \ldots,
  2,1)$, see Lemma~\ref{lem1}. Next use the definition of
  bideterminants $(\di:\jj)$ given in Section
  \ref{definestructuremaps}, Equation (\ref{eqn1}) and sum over all
  $\jj \in I^-(n,r)$ to obtain:
\begin{eqnarray} \label{ast2} 
q^{-|\di|}\sum_{\jj\in I^-(n,r)}(\di:\jj) \cdot \det
=
\sum_{\jj\in I^-(n,r)} q^{-|\jj|} \cdot \det \cdot (\di:\jj)
\end{eqnarray}
for any $\di\in I^-(n,r)$. Let $\tau_{n,r}$ and $\tau_{r,n}$ be the
structure maps of the comodules $\wedge^nV\otimes \wedge^rV$ and
$\wedge^rV\otimes \wedge^nV$ respectively.  Applying these structure
maps to an element $\omega \otimes \hat{v}_{\di}$, we get as
coefficients of an element $\hat{v}_{\underline{j}} \otimes \omega$
precisely the left-hand and right-hand side of Equation (\ref{ast2}).
Then Equation (\ref{ast2}) implies that
\[
((id \otimes \mu_r )\circ \tau_{n,r})(\omega \otimes
\hat{v}_{\di})= (\tau_{r,n} \circ \mu_r)(\omega \otimes
\hat{v}_{\di})
\]
where $\di \in I^-(n,r)$.
So $\phi$ is a $k[\M]$-comodule homomorphism.
\end{pf}
\begin{lem}\label{Mmodulemaps}
  Fix a natural number $m$ and let $\alpha \in P^+(n,r)_m$. Then the linear map $\phi_{\alpha} :
  \wedge^{\alpha} V\otimes \wedge^{n-\alpha}V\To (\wedge^nV)^{\otimes
    m}$ defined by
\[
\phi_\alpha(\hv_S\otimes \hv_{\widehat{T}}) =
\sgn(S)\delta_{S,T}q^{-\epsilon(S)}\omega^{\otimes m}
\]
for any $S, T\in Tab^-(\alpha)$ is an $\M$-module homomorphism.
\end{lem}
\begin{pf}
  To show the claim, we proceed by induction on $m$.  The
  multiplication map $\wedge^s V\otimes \wedge^t V\To \wedge^{s+t}V$
  is an $\M$-module homomorphism. If $m =1$, then $\phi_{\alpha}$ is
  the multiplication map $\wedge^{\alpha_1}V\otimes
  \wedge^{n-\alpha_1}V\To \wedge^nV$ which, by definition, maps
\[
\hv_{\di}\otimes \hv_{\hat{\jj}} \mapsto
\delta_{\di,\jj}\sgn(\di)\omega,
\]
and here $S=\di$ and $T=\jj$ are elements in $Tab^-(\alpha) =
I^-(n,\alpha_1)$.

Now suppose $m>1$.  Given a partition $\alpha=(\alpha_1, \ldots,
\alpha_m) \in P^+(n,r)_m$, assume that for $(\alpha_2, \ldots, \alpha_m)
\in P^+(n,r-\alpha_1)_{m-1}$ the induction assumption holds. Split any
tableau $S \in Tab^-(\alpha)$ into the first column $S_1$ and a
tableau $S_e \in P^+(n,r-\alpha_1)_{m-1}$ with columns $S_2, \ldots,
S_m$. Then $\epsilon(S) = (m-1) |S_1| + \epsilon(S_e)$, and
$\phi_{\alpha}$ is the composite map
\begin{align*}
\wedge^{\alpha}V\otimes \wedge^{n-\alpha}V
&\To \wedge^{\alpha_1}V\otimes (\wedge^{\alpha_2}V\otimes
\ldots\otimes \wedge^{n-\alpha_2}V)\otimes \wedge^{n-\alpha_1}V\\
&\To \wedge^{\alpha_1}V\otimes (\wedge^nV)^{\otimes m-1}\otimes
\wedge^{n-\alpha_1}V \\ &\To \wedge^{\alpha_1}V\otimes
\wedge^{n-\alpha_1}V\otimes (\wedge^{n}V)^{\otimes m-1}\\
&\To (\wedge^{n}V)^{\otimes m}.
\end{align*}
In fact, if $S \neq T$, then both $\phi_{\alpha}$ and the composite
map return zero. Suppose next that $S = T$. By induction and Lemma
\ref{lem3}, the map $\phi_{\alpha}$ is an $\M$-module homomorphism.
\end{pf}

\begin{thm} \label{compareGmaps} Fix a natural number $m$ and
  let $\alpha, \beta$ be partitions in $\W_m$. \\
  There is a vector space isomorphism $\theta_{\alpha}:
  \wedge^{\alpha}V\To \wedge^{n-\alpha}E$ defined by
\[
\theta_{\alpha}(\hv_S) =  \sgn(S)q^{-\epsilon(S)}\he_{\widehat{S}}
\]
for $S \in Tab^-(\alpha)$.
Sending $\varphi$ to
$\theta_{\beta}\circ\varphi\circ\theta^{-1}_{\alpha}$ defines an
isomorphism
\begin{eqnarray} \label{ast}
\Hom_{S(n,r)}(\wedge^{\alpha}V,\wedge^{\beta}V) \simeq
\Hom_{S(n,nm-r)}(\wedge^{n-\alpha}E,\wedge^{n-\beta}E).
\end{eqnarray}
\end{thm}
In other words, the following commutative diagram provides a bijection
between $\M$-morphisms $\varphi$ in degree $r$ and $\M$-morphisms
$\psi$ in degree $nm-r$.
{
\[
\xymatrix{ \wedge^\alpha V \ar[dd]_\varphi \ar[rr]^{\theta_\alpha}
& &\wedge^{n-\alpha}E
\ar[dd]^{\psi:=\theta_{\beta}\circ\varphi\circ\theta^{-1}_{\alpha}} \\
& & \\
\wedge^\beta V \ar[rr]^{\theta_\beta}
& & \wedge^{n-\beta} E
}
\]
}
\begin{pf}
  The map $\theta_{\alpha}$ is an isomorphism of vector spaces; this
  is immediate from the above description of bases for
  $\wedge^{\alpha}V$ and $\wedge^{n-\alpha}E$.  Clearly, the inverse
  of the map $\varphi \mapsto \psi$ is given by $\psi \mapsto
  \theta^{-1}_{\beta}\circ\psi\circ\theta_{\alpha}$. So to prove the
  main statement, that is Equation (\ref{ast}), it suffices to verify
  that $\psi:=\theta_{\beta}\circ\varphi\circ\theta^{-1}_{\alpha}$ is
  an $\M$-module homomorphism if $\varphi$ is so.  As usual, we are
  viewing $\M$-modules as $k[\M]$-comodules.

{\bf Notation:}
Throughout we will be using that indexing over a set of tableaux is
equivalent to indexing over the set of complement tableaux, as defined in
Section~\ref{sec3}.  Suppose that
\[
Tab^-(\alpha) = \{A_1,\ldots,A_a\} \mbox{ and }
Tab^-(\beta)  = \{B_1,\ldots,B_b\}.
\]
Write ${A}$ for the $a\times a$ matrix with $(i,j)^{th}$ entry
$(A_i:A_j)$, and ${B}$ for the $b\times b$ matrix with $(i,j)^{th}$
entry $(B_i: B_j)$.  Moreover, define the following matrices:
\[
\begin{array}{lll}
\sgn(A) &=& diag(\sgn(A_1),\ldots,\sgn(A_a)), \\
\sgn(B) &=& diag(\sgn(B_1),\ldots,\sgn(B_b)),\\
q^A     &=& diag(q^{\epsilon(A_1)},\ldots,q^{\epsilon(A_a)}),  \\
q^B     &=& diag(q^{\epsilon(B_1)},\ldots,q^{\epsilon(B_b)}). \\
\end{array}
\]
Here $\epsilon(-)$ is defined as in Section 3.3.  Since $q \neq
0$, the matrices $q^A$ and $q^B$ have inverses, denoted by $q^{-A}$
and $q^{-B}$ respectively.
Note that
\[
Tab^-(n-\alpha)=\{\widehat{A_1},\ldots,\widehat{A_a}\} \mbox{ and }
Tab^-(n-\beta)=\{\widehat{B_1},\ldots,\widehat{B_b}\}.
\]
Thus we can define matrices $\widehat{A}$,$\widehat{B}$, $\sgn(\widehat{A})$
etc.~for partitions $n - \alpha$ and $n-\beta$ in the same way as we
defined matrices $A$, $B$, $\sgn(A)$ etc.~for $\alpha$ and $\beta$.

{\bf Step 1:} Suppose $\varphi$ is an $\M$-module homomorphism
$\wedge^\alpha V \rightarrow \wedge^\beta V$ given by
\begin{equation} \label{eqnvarphi}
\varphi(\hv_{A_i}) = \sum_{B_j \in Tab^-(\beta)} x_{A_iB_j}\hv_{B_j},
\end{equation}
for $A_i \in Tab^-(\alpha)$. Then by definition, $\psi$ is given by
\begin{equation} \label{eqnpsi}
\psi(\he_{\widehat{A_i}}) = \sum_{B_j\in Tab^-(\beta)}
\sgn(A_i)q^{\epsilon(A_i)}x_{A_iB_j}\sgn(B_j)q^{-\epsilon(B_j)}
\he_{\widehat{B_j}}.
\end{equation}
Let $X$ be the $a\times b$ matrix with $(i,j)^{th}$ entry $x_{A_iB_j}$
given by Equation (\ref{eqnvarphi}).

{\bf Step 2:} By definition, the map $\psi$ is an $\M$-module
homomorphism if and only if $(\tau_{n-\beta} \circ \psi)=((\psi
\otimes 1) \circ \tau_{n-\alpha})$; here $\tau_{n-\alpha}$ and
$\tau_{n-\beta}$ denote the comodule structure maps of
$\wedge^{n-\alpha}E$ and $\wedge^{n-\beta}E$ respectively, see Section
2.3.  In the following we spell out this latter equation.  Let
$\hat{e}_{\widehat{A_i}}$ be a basis element in $\wedge^{n-\alpha}E$.  By
Equation (\ref{eqnpsi}) and Lemma~\ref{lem2}, we get
\begin{eqnarray*}
\sum_{j=1}^b  \sum_{l=1}^b
&& \sgn(A_i) q^{\epsilon(A_i)}x_{A_iB_l}\sgn(B_l) q^{-\epsilon(B_l)}
\hat{e}_{\widehat{B_j}} \otimes
\langle \widehat{B_j}:\widehat{B_l} \rangle \\
&=& (\tau_{n-\beta} \circ \psi) (\hat{e}_{\widehat{A_i}})\\
&=&((\psi \otimes 1) \circ \tau_{n-\alpha}) (\hat{e}_{\widehat{A_i}})\\
&=&\sum_{j=1}^b \sum_{u=1}^a
\sgn(A_u)q^{\epsilon(A_u)}x_{A_uB_j}\sgn(B_j)q^{-\epsilon(B_j)}
\hat{e}_{\widehat B_j} \otimes
\langle \widehat{A_u}:\widehat{A_i} \rangle.
\end{eqnarray*}
This equation holds precisely if equality holds for the
coefficients of each basis element $\he_{\widehat B_j}$. Hence the map
$\psi$ is an $\M$-module homomorphism if and only if the following
matrix identity holds:
\begin{eqnarray*} \label{matrixeqntoprove}
& & \sgn(A) q^A X \sgn(B) q^{-B}\widehat{B}^{tr} \\&=& \left(
\sum_{l=1}^b \sgn(A_i) q^{\epsilon(A_i)}x_{A_iB_l}\sgn(B_l)
q^{-\epsilon(B_l)} \langle
\widehat{B_j}:\widehat{B_l} \rangle \right)_{i=1 \ldots a,\,\, j=1 \ldots b}\\
&=& \left(\sum_{u=1}^a
\sgn(A_u)q^{\epsilon(A_u)}x_{A_uB_j}\sgn(B_j)q^{-\epsilon(B_j)} \langle
\widehat{A_u}:\widehat{A_i} \rangle \right)_{i=1 \ldots a,\,\, j=1 \ldots b} \\
&=& \widehat{A}^{tr} \sgn(A) q^A X \sgn(B) q^{-B}
\end{eqnarray*}
Note that the matrix $X$ in the middle of the expressions determines
what are rows and what are columns; thus the matrices $\widehat{A}$ and
$\widehat{B}$ have to be transposed. Multiplying this equation from the
left by $\sgn(A)$ and on the right by $\sgn(B)$ gives the equivalent
equation
\begin{eqnarray} \label{eqn6}
\Pi:=
q^A X \sgn(B) q^{-B}\widehat{B}^{tr} \sgn(B)
= \sgn(A) \widehat{A}^{tr} \sgn(A) q^A X q^{-B}
=:\Gamma.
\end{eqnarray}
Hence the map $\psi$ is an $M$-module homomorphism if and only if
$\Pi=\Gamma$.

{\bf Step 3:} The map $\varphi$ given by Equation (\ref{eqnvarphi}) is
an $\M$-module homomorphism by assumption. Hence, by definition,
$\tau_\beta \circ \varphi=(id \otimes \varphi)\circ \tau_\alpha$.  Let
$\hv_{A_i}$ be a basis element of $\wedge^\alpha V$. By
Lemma~\ref{lem2} and Equation (\ref{eqnvarphi}) we have
\begin{eqnarray*}
\sum_{t=1}^b\sum_{j=1}^b x_{A_iB_t} (B_t:B_j)\otimes \hv_{B_j} &=&
(\tau_\beta \circ \varphi)(\hv_{A_i}) \\ &=& ((id \otimes
\varphi)\circ \tau_\alpha) (\hv_{A_i}) \\ &=&\sum_{u=1}^a
\sum_{j=1}^b  ({A_i}:{A_u}) x_{A_uB_j}  \otimes \hv_{B_j}.
\end{eqnarray*}
Comparing the coefficients of the basis elements $\hv_{B_j}$, we get
that the following matrix equation is satisfied:
\begin{eqnarray}  \label{eqn3}
  XB = \left( \sum_{t=1}^b x_{A_iB_t} (B_t:B_j)
\right)_{i,j}
= \left(\sum_{u=1}^a (A_i:A_u) x_{A_uB_j}\right)_{i,j}
=AX.
\end{eqnarray}
{\bf Step 4:} On $(\wedge^nV)^{\otimes m}$, the structure map
$\tau_n^{\otimes m}=\tau_{(n, \ldots, n)}$ is multiplication by
$(det)^m$, see Lemma \ref{lem2}(3).  We use the map
$\phi_{\alpha}$ defined in Lemma \ref{Mmodulemaps}. Let $\hv_{A_i}
\otimes \hv_{\widehat{A_j}}$ be a basis element in $\wedge^\alpha V
\otimes \wedge^{n-\alpha}V$. Then, since $\phi_{\alpha}$ is an
$\M$-module homomorphism, and using Lemma \ref{lem2}:
\begin{eqnarray*} 
& & \sum_{s=1}^a (A_i:A_s) (\widehat{A_j}:\widehat{A_s})
 \sgn(A_s) q^{-\epsilon(A_s)} \omega^{\otimes m}\\
&=& \sum_{s=1}^a \sum_{t=1}^a (A_i:A_s) (\widehat{A_j}:\widehat{A_t})
 \delta_{A_s,A_t} \sgn(A_s) q^{-\epsilon(A_s)} \omega^{\otimes m} \\
 &=&((id \otimes \phi_{\alpha}) \circ (\tau_{\alpha} \otimes \tau_{n- \alpha}))
 (\hv_{A_i} \otimes \hv_{\widehat{A_j}}) \\ &=& (\tau_{n}^{\otimes m}
 \circ \phi_{\alpha})(\hv_{A_i} \otimes \hv_{\widehat{A_j}})\\ &=&
\delta_{A_i,A_j} \sgn(A_i) q^{-\epsilon (A_i)} (det)^m \omega^{\otimes m}.
\end{eqnarray*}
Comparing the coefficients of $\omega^{\otimes m}$ in this last
equation, we obtain the matrix equation:
\begin{eqnarray*}
A \sgn(A) q^{-A} \widehat{A}^{tr} & = & \left( \sum_{s=1}^a (A_i:A_s)
(\widehat{A_j}:\widehat{A_s}) \sgn(A_s) q^{-\epsilon(A_s)}\right)_{i,j} \\
 &=&\left( \delta_{A_i,A_j} \sgn(A_i) q^{-\epsilon (A_i)} (det)^m
\right)_{i,j} \\
&=& \sgn(A) q^{-A}(\det)^m I.
\end{eqnarray*}
Since the quantum parameter $q$ is non-zero, the matrix on the
right-hand side of the last equation is invertible in
$k[\M][{\det}^{-1}]$.  This implies that the matrix $A$ is right
invertible and the matrix $\widehat{A}$ is left invertible in
$k[\M][{\det}^{-1}]$. Reversing the roles of $A$ and $\widehat{A}$ and
using that $\widehat{\widehat{A}} =A$, we get that $A$ is invertible on either
side in $k[\M][{\det}^{-1}]$.  We next multiply this equation from the
left by $q^{A}$ and from the right by $\sgn(A)$, obtaining:
\begin{eqnarray}
q^A A \sgn(A) q^{-A} \widehat{A}^{tr}\sgn(A) &=& (\det)^m I, \label{eqn4}\\
q^B B \sgn(B)q^{-B} \widehat{B}^{tr}\sgn(B) &=& (\det)^m I,  \label{eqn5}
\end{eqnarray}
where the second equation is obtained similarly by using the map
$\phi_{\beta}$ instead of $\phi_{\alpha}$.

{\bf Step 5:} Our aim is to prove that $\psi$ is an $M$-module
homomorphism.  By Step 2 it is enough to show that Equation
(\ref{eqn6}) holds, that is $\Gamma=\Pi$. In order to show the
latter, we next combine the matrix equalities obtained in previous
steps:
\begin{align*}
q^A A q^{-A} \Gamma
&= [q^A A q^{-A} \sgn(A) \widehat{A}^{tr} \sgn(A) ] q^A X q^{-B} \\
&= (\det)^m q^A X q^{-B} \mbox{\hspace{5cm} (using Equation (\ref{eqn4}))} \\
&= q^A X q^{-B} (\det)^m \\
&= q^A X q^{-B} q^B B q^{-B} \sgn(B) \widehat{B}^{tr} \sgn(B)
\mbox{\hspace{1.7cm} (using Equation (\ref{eqn5}))} \\
&= q^A A q^{-A} [q^A X \sgn(B) q^{-B} \widehat{B}^{tr} \sgn(B) ]
\mbox{\hspace{1.5cm} (using Equation (\ref{eqn3}))} \\
&= q^A A q^{-A} \Pi.
\end{align*}
Hence $0 = q^{A}Aq^{-A}(\Gamma - \Pi)$, and since $A$ is invertible by
Step 4, this implies the validity of Equation (\ref{matrixeqntoprove})
in $k[\M][{\det}^{-1}]$, and hence in the subcoalgebra $k[\M]$.
\end{pf}
An obvious, but crucial, corollary of Theorem~\ref{compareGmaps}
provides isomorphisms of rings.  Here, as always we assume $q\neq
0$; hence we do not have to distinguish between $V$ and $E$, which are
related by an $\M$-pairing, see \cite{Donkin}.
\begin{cor} \label{ringiso}
For any subset $X\subset \W_m$, there is an isomorphism of algebras:
\[
\End_{S(n,r)}(\bigoplus_{\alpha \in X} \wedge^{\alpha} V) \simeq
\End_{S(n,nm-r)}(\bigoplus_{\alpha \in X} \wedge^{n-\alpha}V).
\]
\end{cor}

\noindent{\bf Remark:}
This corollary implies equivalences of categories, as stated in the
introduction. These equivalences may be obtained more directly. The
category of all representations of the (algebraic or quantised) group
G has an autoequivalence induced by inverting group elements (that is,
by the Hopf algebra antipode) and transposing matrices. Other
autoequivalences can be obtained by tensoring with powers of the
determinant. Appropriately composing these two kinds of equivalences
and restricting to homogeneous polynomial representations produces the
equivalences implied by the Corollary.
Adding combinatorial arguments, as developed in Section 3, would then
prove the existence of isomorphisms as in the previous Corollary,
without explicitly constructing them.

We also note that the set $X$ in the Corollary can be chosen arbitrarily.
We will only use a special choice of $X$, which is compatible with the
quasi-hereditary structure.

\section{Ringel Duality}

\subsection{General results on Ringel duality}
Tensor products of exterior powers, the representations used in the
previous section, have a special structural property: they are direct
summands of the characteristic tilting module of the Schur algebra in
question. The present section uses this property, and known results on
tilting modules, to turn Corollary \ref{ringiso} into a statement on
Schur algebras.  Let us recall some basics on quasi-hereditary
algebras and Ringel duality. A convenient reference in our context is
the appendix of \cite{Donkin}.

\begin{defn}
Let $\{L(\lambda) :\lambda\in \Lambda\}$ be a complete set of
simple modules of a finite dimensional algebra $A$, $P(\lambda)$
the projective cover of $L(\lambda)$ and let $(\Lambda, \leq)$ be
a partial ordering on $\Lambda$. Then $(A, \Lambda,\leq)$ is
called a {\em quasi-hereditary algebra} if for each $\lambda\in
\Lambda$, there exists a quotient module $\Delta(\lambda)$ of
$P(\lambda)$, called {\em standard module}, such that
\begin{enumerate}
    \item  The kernel of the canonical map $P(\lambda)\twoheadrightarrow
    \Delta(\lambda)$ is filtered by $\Delta(\mu)$ with $\mu>
    \lambda$;
    \item The kernel of the canonical map $\Delta(\lambda)\twoheadrightarrow
    L(\lambda)$ is filtered by $L(\mu)$ with $\mu< \lambda$.
\end{enumerate}
\end{defn}

An equivalent ring-theoretic definition for quasi-hereditary algebras
implies that $A$ is quasi-hereditary if and only if $A^{op}$ is so. If
we denote the standard module of $A^{op}$ by
$\Delta_{A^{op}}(\lambda)$, then $\nabla(\lambda) = \nabla_A(\lambda)
= \Delta_{A^{op}}(\lambda)^*$ is called a {\em costandard module} of
$(A,\leq)$. Here $^*$ denotes the usual vector space dual.

Let $\mathcal{F}(\Delta)$ be the full subcategory of $A$-mod
consisting of modules filtered by $\Delta$'s and let
$\mathcal{F}(\nabla)$ be the full subcategory of $A$-mod consisting of
modules filtered by $\nabla$'s. Ringel \cite{R} proved that for each
$\lambda\in \Lambda$ there is an indecomposable module $T(\lambda)$
such that
\[
\mathcal{F}(\Delta)\cap \mathcal{F}(\nabla) = add(\oplus_{\lambda\in
\Lambda} T(\lambda)).
\]
The module $T = \oplus_{\lambda\in \Lambda}T(\lambda)$ is called the
{\em characteristic tilting module} of $(A,\leq)$. The endomorphism
ring $R = \End_A(T)$, and any algebra Morita equivalent to it, is
called {\em Ringel dual} of $A$.  Ringel has shown that $R$ is again a
quasi-hereditary algebra.  Moreover, the Ringel dual of $R$ is Morita
equivalent to $A$ (see \cite{R}).

We need one more definition in this context: A subset $Y$ of $\Lambda$
is called {\it saturated} if $y_2<y_1$ and $y_1\in Y$ implies $y_2\in
Y$.  The following statements are well-known and easy to prove.
\smallskip

\begin{lem} \label{Ringelduality}
  Let $(A,\Lambda,\leq)$ be a quasi-hereditary algebra and let $T =
  \oplus_{\lambda \in\Lambda } T(\lambda)$ be its characteristic
  tilting module. For any saturated subset $Y$ of $\Lambda$,
\begin{enumerate}
    \item the primitive idempotents $e_{\lambda}$ with $\lambda\notin
    Y$ generate an ideal $J$ in the heredity chain of $A$;
    \item the module $T(Y) = \oplus_{\lambda\in Y} T(\lambda)$ is
    exactly the characteristic tilting module of the quasi-hereditary
    algebra $A/J$;
    \item the endomorphism ring $\End_A(T(Y)) = \End_{A/J}(T(Y))$
    is of the form $eRe$, where $R$ is the Ringel dual of $A$ and
    $e$ is an idempotent in $R$.
\end{enumerate}
\end{lem}

The (classical or quantized) Schur algebras $S(n,r)$ are
quasi-hereditary algebras with respect to $\Lambda=\Lambda^+(n,r)$ and
the dominance ordering $\leq$.  

\begin{thm}[Donkin \cite{D}, \cite{Donkin}, Section 3.3]
\label{Donkinstheorem}
The indecomposable tilting modules for the Schur algebra $S(n,r)$ are
precisely the direct summands of the modules $\wedge^{\alpha'} E$, for
$\alpha\in \LW(n,r)$. The indecomposable tilting module $T(\alpha)$ occurs
exactly once as a direct summand of $\wedge^{\alpha'} E$, and if
$T(\lambda)$ is a direct summand of $\wedge^{\alpha'} E$ then $\lambda
\leq \alpha$.
\end{thm}

\subsection{Morita equivalences between generalised Schur algebras}
Given a Schur algebra $S(n,r)$ and a saturated set $\Pi\subset
\Lambda^+(n,r)$.  Define $I_{\Pi}$ to be the ideal in the heredity
chain of $S(n,r)$, generated by the primitive idempotents $e_\alpha$
with $\alpha \notin \Pi$. Moreover, after fixing a natural number $m$
such that the partitions in $\Pi$ have rows of length not more than
$m$, define $\widehat{\Pi}$ as the set with elements
$\hat{\alpha}:=T_{m}(\alpha)$ for $\alpha \in \Pi$. Note that
$\widehat{\Pi}$ is again saturated.
Combining Corollary \ref{ringiso} with Lemma \ref{Ringelduality}
and Theorem \ref{Donkinstheorem} we can show:
\begin{thm} \label{RingeldualMorita}
  Fix a natural number $m$.  For any natural numbers $n$ and $r$ and
  any saturated set $\Pi\subset \Lambda^+(n,r)$, such that the
  partitions in $\Pi$ have rows of length not more than $m$, there is
  a Morita equivalence:
\[
S(n,r)/I_{\Pi} \simeq S(n,nm-r)/I_{\widehat{\Pi}}.
\]
\end{thm}
\begin{pf}
  Consider the following modules over $S(n,r)$ and $S(n, nm-r)$
  respectively:
\begin{align*}
T
&=\bigoplus_{\alpha\in \Pi}\wedge^{\alpha'}
E \oplus \bigoplus_{\beta\notin \Pi}T(\beta),\\
\tilde{T}
&=\bigoplus_{\alpha\in\Pi}\wedge^{(n-\alpha')}
E\oplus \bigoplus_{\beta\notin\widehat{\Pi}}T(\beta)
\end{align*}
which are full tilting modules by Theorem \ref{Donkinstheorem}.  Note
that the exterior powers here are defined using conjugate partitions
$\alpha'$ for $\alpha \in \Pi$.  Let
\[
A =\End_{S(n,r)}(T) \mbox{ and }
B = \End_{S(n,nm-r)}(\tilde{T})
\]
be the Ringel duals of $S(n,r)$ and $S(n,nm-r)$ respectively.  By
Corollary \ref{ringiso}, we have
\[
eAe = \End_{S(n,r)}(\oplus_{\alpha\in \Pi}\wedge^{\alpha'} E)\simeq
\End_{S(n,nm-r)}(\oplus_{\alpha\in \Pi}\wedge^{(n-\alpha')}E) = fBf
\]
where $e$ and $f$ are idempotents in $A$ and $B$. By Lemma
\ref{Ringelduality} and Theorem \ref{Donkinstheorem}, the direct sum
$\oplus_{\alpha\in \Pi}\wedge^{\alpha'} E$ is a full tilting module
over $S(n,r)/I_{\Pi}$ with $I_{\Pi}$ an ideal in the heredity chain of
$S(n,r)$, generated by the primitive idempotents $e_\alpha$ with
$\alpha\notin \Pi$. Moreover $\oplus_{\alpha\in \Pi}\wedge^{n-\alpha'}
E$ is a full tilting module over $S(n,nm-r)/I_{\widehat{\Pi}}$ with
$I_{\widehat{\Pi}}$ an ideal in the heredity chain of $S(n, nm-r)$,
generated by the primitive idempotents $e_\beta$ with $\beta\notin
\widehat{\Pi}$. Hence, using Ringel duality, we get the following
Morita equivalence:
\[
S(n,r)/I_{\Pi} \simeq S(n,nm-r)/I_{\widehat{\Pi}}.
\]
\end{pf}

\section{Isomorphisms between generalized Schur algebras}
It just remains to be shown that the Morita equivalences in Theorem
\ref{RingeldualMorita} indeed are isomorphisms between generalized
Schur algebras.

\begin{lem}\label{equaldimiso}
  Let $(A,\leq_A)$ and $(B,\leq_B)$ be two quasi-hereditary $k$-algebras,
  which are Morita equivalent such that the Morita equivalence induces
  an isomorphism between $\leq_A$ and $\leq_B$. Then this Morita
  equivalence sends standard modules to standard modules, say
  $\Delta_A(i)$ to $\Delta_B(i)$. Furthermore, suppose for all $i$,
  the modules $\Delta_A(i)$ and $\Delta_B(i)$ have the same
  $k$-dimension. Then the algebras $A$ and $B$ are isomorphic.
\end{lem}
\begin{pf} Standard modules are relatively projective, namely
projective in truncated categories, that is as modules over
quasi-hereditary quotients. Hence the Morita equivalence preserves the
property of being standard.
The Morita equivalence also sends simple modules to simple modules,
thus it identifies the decomposition matrices $([\Delta(i):
L(j)])_{i,j}$ of $A$ and of $B$. By induction it follows that the
$A$-simple $L_A(i)$ has the same $k$-dimension as the $B$-simple
$S_B(i)$.
Hence the Morita equivalence does not change multiplicities of
projectives in the regular representation and therefore it
provides an algebra isomorphism.
\end{pf}

We keep the notation of the previous section for classical or
quantized Schur algebras. Let $\alpha \in \Lambda^+(n,r)$ with
$\alpha_1 \leq m$.

\begin{lem} \label{mirrodimensions}
  The two Weyl modules $\Delta( \alpha)$ and $\Delta(\hat{\alpha})$
  have the same $k$-dimension.
\end{lem}
\begin{pf}
  It is known (see \cite{Green, Donkin}) that the $k$-dimension of the
  Weyl module $\Delta(\alpha)$ is precisely the number of
  semi-standard $\alpha$-tableaux: $\dim_k\Delta(\alpha) =
  |STab(\alpha)|$.  Note that in the literature semistandard tableaux
  are often defined by requiring entries to increase along rows and to
  increase strictly along columns, whereas we have required entries to
  (weakly/strictly) decrease.  Replacing an entry $i$ by $n+1-i$ provides a
  bijection between the two kinds of semistandard tableaux. By Theorem
  \ref{ssmirror} we have $|STab(\alpha)| = |STab(\hat{\alpha})|$.
\end{pf}

An alternative proof of this Lemma could be based on the isomorphism
$\Delta( \hat{\alpha}) \simeq \nabla(\alpha)^{\ast} \otimes det^{\otimes
  m}$. This can be shown by checking that both sides have the same highest
weight and also the same dimension, by Weyl's character formula.
\bigskip

Recall the definition of the ideals $I_{\Pi}$ and $I_{\widehat{\Pi}}$
given in Theorem~\ref{RingeldualMorita}.  We get the main
theorem:

\begin{thm} \label{comparisontheorem} \label{maintheorem}
  For any natural numbers $n$ and $r$ and for any saturated set
  $\Pi\subset \Lambda^+(n,r)$,
  and for any natural number $m$ such that the partitions in $\Pi$
  have rows of length not more than $m$, there is an isomorphism
  between (classical or quantized) generalized Schur algebras:
$$S(n,r)/I_{\Pi} \simeq S(n,nm-r)/I_{\widehat{\Pi}}.$$
The resulting equivalence of categories of polynomial
$G$-modules sends $M(\alpha)$ to $M(\hat{\alpha})$ where $M \in
\{L,\Delta, P,T \}$.
\end{thm}
\begin{pf}
  This is a consequence of Theorem \ref{RingeldualMorita}, Lemma
  \ref{equaldimiso} and Lemma \ref{mirrodimensions}.
\end{pf}

We note that the isomorphisms in Theorem \ref{comparisontheorem} in
the classical case do not depend on the choice of the ground field or
on its characteristic, and in the quantized case they do not depend on
the value of the quantum parameter $q$, provided $q \neq 0$.

The embedding of $S(n,r)$-mod into $\M$-mod preserves cohomology, and
so does for any ideal $I$ in a heredity chain the embedding
$S(n,r)/I$-mod into $S(n,r)$-mod. Therefore Theorem \ref{maintheorem}
provides us with isomorphisms of cohomology groups; using the above notation
we get for example:
\[
Ext^{\ast}(L(\alpha),L(\beta)) \simeq
Ext^{\ast}(L(\hat{\alpha}),L(\hat{\beta})).
\]
Here, the extension groups may be taken over $\M$ or over the respective
(classical or quantized) Schur algebras.

\section{Applications: Decomposition numbers and structure of 
Schur algebras}

We list different applications of the isomorphisms and equivalences of
categories provided by Theorem \ref{maintheorem}: identities of
decomposition numbers, isomorphisms of cohomology groups, 
and a factorisation of the isomorphisms
provided earlier by \cite{BLM}.  Unexplained terminology, background
information and references can be found, for example, in \cite{K}. We
keep the notation of the previous sections. 
\subsection{Decomposition numbers of Schur algebras}
\begin{cor}
The following identity holds for the decomposition numbers
of classical and quantum $GL_n$: Fix a natural number $m$. Then
$$[\Delta(\lambda):L(\mu)] = [\Delta(\hat{\lambda}):L(\hat{\mu})]$$
for any positive integer $r$ and for any $\lambda$ and $\mu$ with
$\lambda',\mu'\in \PW_m$.
\end{cor}
Note that by varying $m$ and iterating the process of taking
complements, we can get many (in general, infinitely many)
different complement partitions $\hat{\lambda}$ and $\hat{\mu}$,
thus equating the decomposition number $[\Delta(\lambda):L(\mu)]$
with many other decomposition numbers in various degrees.

{\bf Example.}  Denote by $D_p(n,r)$ the decomposition matrix of the
classical Schur algebra $S(n,r)$ over a field of prime
characteristic. Label the columns and rows of these decomposition
matrices in the same order, in both cases starting in the top left corner.
Moreover, denote the decomposition number zero by a dot.  Choose $p=2$
and $p=5$ respectively, and consider the Schur algebras $S(3,5)$ and
$S(3,10)$.
\[
\begin{array}{cccc}  
\begin{array}{|c|ccccc|}
\hline
D_2(3,5) &&&&&\\
\hline
( 5, 0, 0 ) &1&.&1&1&.\\
( 4, 1, 0 ) &.&1&.&.&.\\
( 3, 2, 0 ) &.&.&1&1&1\\
( 3, 1, 1 ) &.&.&.&1&1\\
( 2, 2, 1 ) &.&.&.&.&1\\
\hline
\end{array}
&
\mbox{ and }
&
\begin{array}{|c|ccccc|}
\hline
D_5(3,5) &&&&& \\
\hline
( 5, 0, 0 ) & 1&1&.&.&.\\
( 4, 1, 0 ) & .&1&.&1&.\\
( 3, 2, 0 ) & .&.&1&.&.\\
( 3, 1, 1 ) & .&.&.&1&.\\
( 2, 2, 1 ) & .&.&.&.&1\\
\hline
\end{array}
\end{array}
\]
%
We notice that the decomposition matrix $D_2(3,5)$ occurs as a
submatrix of $D_2(3,10)$; in the latter matrix, we consider the
submatrix defined by the rows and corresponding columns marked with
a star. Similarly for prime $5$, the decomposition matrices $D_5(3,5)$
is contained in $D_5(3,10)$ in the same way as for $p=2$.
{\footnotesize
\[
\begin{array}{|l|lllll|l|lll|l|l|lll|} 
\hline D_2(3,10)&&&&&&*&&&&*&&*&*&* \\ \hline
(10, 0, 0) &1&1&.&1&1&1&.&.&.&.&1&.&.&.\\
( 9, 1, 0 ) &.&1&1&1&.&1&1&.&.&.&.&.&.&1\\
( 8, 2, 0 ) &.&.&1&1&.&1&1&.&.&.&1&1&1&1\\
( 7, 3, 0 ) &.&.&.&1&1&1&.&.&1&.&2&1&2&.\\
( 6, 4, 0 ) &.&.&.&.&1&1&.&.&1&.&1&1&2&.\\
\hline
( 5, 5, 0 )^* &.&.&.&.&.&{\bf 1}&.&.&.&{\bf .}&.&{\bf 1}&{\bf 1}&\bf .\\
\hline
( 8, 1, 1 ) &.&.&.&.&.&.&1&.&.&.&1&.&1&1\\
( 7, 2, 1 ) &.&.&.&.&.&.&.&1&.&1&.&.&.&.\\
( 6, 3, 1 ) &.&.&.&.&.&.&.&.&1&.&1&1&2&.\\
\hline
( 5, 4, 1 )^* &.&.&.&.&.&\bf .&.&.&.&\bf 1&.&\bf .&\bf .&\bf .\\
\hline
( 6, 2, 2 ) &.&.&.&.&.&.&.&.&.&.&1&1&1&.\\
\hline
( 5, 3, 2 )^* &.&.&.&.&.&\bf.&.&.&.&\bf .&.&\bf 1&\bf 1&\bf 1\\
( 4, 4, 2 )^* &.&.&.&.&.&\bf .&.&.&.&\bf .&.&\bf .&\bf 1&\bf 1\\
( 4, 3, 3 )^* &.&.&.&.&.&\bf .&.&.&.&\bf .&.&\bf .&\bf .&\bf 1\\
\hline
\end{array}
\]
}

{\footnotesize
\[
\begin{array}{|l|lllll|l|lll|l|l|lll|} 
\hline D_5(3,10)&&&&&&*&&&&*&&*&*&* \\
\hline
(10, 0,0)&1&1&.&.&.&.&.&.&.&.&.&.&1&.\\
(9, 1, 0)&.&1&.&.&.&1&1&.&.&1&.&.&1&.\\
(8, 2, 0)&
.&.&1&.&1&.&.&.&.&.&.&.&.&.\\
(7, 3, 0)&
.&.&.&1&.&.&.&.&.&.&.&.&.&.\\
(6, 4, 0)&
.&.&.&.&1&.&.&.&.&.&.&.&.&.\\
\hline (5, 5, 0)^*&
.&.&.&.&.&{\bf 1}&.&.&.&{\bf 1}&.&.&.&.\\
\hline (8, 1, 1)&
.&.&.&.&.&.&1&.&.&1&.&.&.&.\\
(7, 2, 1)&
.&.&.&.&.&.&.&1&1&.&.&.&.&.\\
(6, 3, 1)&
.&.&.&.&.&.&.&.&1&.&.&.&.&1\\
\hline (5, 4, 1)^*&
.&.&.&.&.&.&.&.&.&{\bf 1}&.&.&{\bf 1}&.\\
\hline (6, 2, 2)&
.&.&.&.&.&.&.&.&.&.&1&.&.&.\\
\hline (5, 3, 2)^*&
.&.&.&.&.&.&.&.&.&.&.&{\bf 1}&.&.\\
(4, 4, 2)^*&
.&.&.&.&.&.&.&.&.&.&.&.&{\bf 1}&.\\
(4, 3, 3)^*&
.&.&.&.&.&.&.&.&.&.&.&.&.&{\bf 1}\\
\hline
\end{array}
\]
}

\noindent
In fact, the Schur algebra $S(3,5)$ is isomorphic to a quotient of the
Schur algebra $S(3,10)$, and this is true over any (infinite) field of
any characteristic: Choose $m=5$ and $\Pi=\Lambda^+(3,5)$. Then
\[
\widehat{\Pi}:=\{ (5,5,0), (5,4,1), (5,3,2), (4,4,2), (4,3,3) \} \subseteq
\Lambda^+(3,10),
\]
which is precisely the set of partitions indexed by a star in the
above decomposition matrices.
\bigskip

Finally we note that on a computational level, Weyl modules associated
to a partition or to its complement have already been compared by
Pittaluga and Strickland \cite{PS}: explicit bases of these Weyl
modules have been compared, in order to get improved bounds for
dimensions of simple modules. Theorem \ref{comparisontheorem} explains
on a structural level why this comparison of bases is natural.

\subsection{Factorisation of maps constructed by
Beilinson, Lusztig and MacPherson} \label{sec4}

In their geometric study of quantum groups of type $A$, Beilinson,
Lusztig and MacPherson \cite{BLM} constructed surjective maps
\[
S(n,r+n) \twoheadrightarrow S(n,r),
\]
whose kernels $J$ have been shown to be ideals in the heredity chain
of the Schur algebra $S(n,r+n)$.  Hence these maps can be viewed as
isomorphisms between generalized Schur algebras: $S(n,r+n)/J \simeq
S(n,r)$.  For more information on these maps see also \cite{Du,
  RMGreen, GR}.

Indeed, these maps can be written as products of two isomorphisms
provided by Theorem \ref{maintheorem} as follows: Start with
$\Lambda^+(n,r+n)$ and $m=r+1$. Let $\lambda$ be a hook partition with
first row containing $m=r+1$ boxes and first column containing $n$
boxes. Define $\Pi$ to be the saturated set of all partitions $\mu$
less than or equal to $\lambda$ in the dominance ordering. These are
precisely the partitions in $\Lambda^+(n,r+n)$ whose first column has
exactly $n$ boxes.  Applying the complement construction to $\lambda$
yields a partition $\hat{\lambda}=n-\lambda \in \Lambda^+(n,(n-1)r)$,
which has the form of an $n \times r$ rectangle. The partitions in
$\widehat{\Pi}$ are precisely those partitions $\hat{\mu}$ of $(n-1)r$
into not more than $n$ parts, which are less than or equal to
$\hat{\lambda}$, since the complement construction preserves the
dominance order, by Lemma \ref{prop1}.  Applying the complement
construction for a second time we now use $m'=r=m-1$. This sends
$\hat{\lambda}$ to $\hat{\hat{\lambda}}$, which has just one row with
$r$ boxes. The partitions in $\widehat{\Pi}$ are being sent to the
partitions $\hat{\hat{\mu}}$ of $r$ into not more than $n$ parts,
which are less than or equal to $\hat{\hat{\lambda}}$; this set is all
of $\Lambda^+(n,r)$.

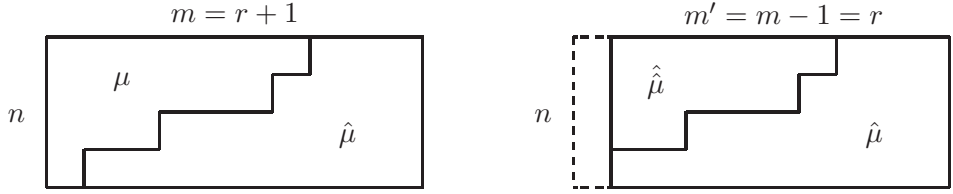
$$
\vspace{1.7cm}
\begin{minipage}{132mm}
  \setlength{\unitlength}{1mm}

\begin{picture}(0,0)(-18,17)
\thicklines

\drawline[20](-10,20)(40,20)
\drawline[20](-10,0)(40,0)
\drawline[20](-10,20)(-10,0)
\drawline[20](65,20)(65,0)
\drawline[20](40,0)(40,20)
\drawline[20](110,0)(110,20)
\drawline[20](65,0)(110,0)
\drawline[20](65,20)(110,20)

\drawline[20](60,0)(60.5,0)
\drawline[20](60,20)(60.5,20)
\drawline[20](61,0)(62,0)
\drawline[20](61,20)(62,20)
\drawline[20](63,0)(64,0)
\drawline[20](63,20)(64,20)
\drawline[20](64.5,0)(65,0)
\drawline[20](64.5,20)(65,20)

\drawline[20](60,20)(60,19)
\drawline[20](60,18)(60,17)
\drawline[20](60,16)(60,15)
\drawline[20](60,14)(60,13)
\drawline[20](60,12)(60,11)
\drawline[20](60,10)(60,9)
\drawline[20](60,8)(60,7)
\drawline[20](60,6)(60,5)
\drawline[20](60,4)(60,3)
\drawline[20](60,2)(60,1)

\drawline[20](-5,0)(-5,5)
\drawline[20](-5,5)(5,5)
\drawline[20](5,5)(5,10)
\drawline[20](5,10)(20,10)
\drawline[20](20,10)(20,15)
\drawline[20](20,15)(25,15)
\drawline[20](25,15)(25,20)

\drawline[20](65,0)(65,5)
\drawline[20](65,5)(75,5)
\drawline[20](75,5)(75,10)
\drawline[20](75,10)(90,10)
\drawline[20](90,10)(90,15)
\drawline[20](90,15)(95,15)
\drawline[20](95,15)(95,20)

\put(-14,9.5){\makebox(0,0){$n$}}
\put(56,9.5){\makebox(0,0){$n$}}
\put(15,23){\makebox(0,0){$m=r+1$}}
\put(88,23){\makebox(0,0){$m'=m-1=r$}}
\put(0,14){\makebox(0,0){$\mu$}}
\put(30,7){\makebox(0,0){$\hat{\mu}$}}
\put(71,14){\makebox(0,0){$\hat{\hat{\mu}}$}}
\put(100,7){\makebox(0,0){$\hat{\mu}$}}

\end{picture}
\end{minipage}
$$

Altogether we have factored the isomorphism from \cite{BLM} as follows:
\[
S(n,r+n)/J \simeq S(n,(n-1)r)/J' \simeq S(n,r)
\]
where $J$ and $J'$ are ideals in the respective heredity chains. So in
particular, James' column removal formula is recovered:

{\bf Column Removal:} For a partition $\lambda=(\lambda_1, \lambda_2,
  \ldots, \lambda_n)$ with non-zero $\lambda_n$, define
  $\tilde{\lambda}=(\lambda_1-1, \ldots, \lambda_n-1)$.  Let
  $\lambda$, $\mu$ be partitions of $r$ such that $\lambda_n$ and
  $\mu_n$ are non-zero.  Then
  $$[\Delta(\lambda) : L(\mu)] = [\Delta({\tilde{\lambda}})
  :L({\tilde{\mu}})].$$

Moreover, iterating this construction will relate Schur algebras
$S(n,r)$ with generalized Schur algebras $S(n,r+ln)/J$, for all
natural numbers $l$ (see \cite{BLM}).
%

\section{Application: Row removal in Schur algebras} \label{sec5}
In this section, we reprove James' row removal formula for
decomposition numbers (see \cite{james}) as an application of
Theorem~\ref{maintheorem} combined with Green \cite{Green}, Section
6.5, showing in particular, that this is not just a numerical result
but reflects an isomorphism between generalized Schur algebras.

{\bf Row Removal:}
  For a partition $\lambda=(\lambda_1, \lambda_2, \ldots, \lambda_n)$,
  define $\tilde{\lambda}=(\lambda_2, \ldots, \lambda_n)$.  Let
  $\lambda$, $\mu$ be partitions of $r$ such that $\lambda_1=\mu_1$.
  Then
  $$[\Delta(\lambda) : L(\mu)] = [\Delta({\tilde{\lambda}})
  :L({\tilde{\mu}})].$$
  Before constructing the underlying isomorphism of generalized Schur
  algebras which in particular also equates decomposition numbers, we
  recall how to obtain row removal numerically, using the above
  strategy. Compare with Donkin's proof of James's row removal
  (see~\cite{donkin1}); also note that row removal is just a special
  case of Donkin's horizontal cut principal, given in~\cite{donkin2},
  Theorem 1 and~\cite{Donkin}, Equation 4.2(9).
\subsection{Row removal numerically} \label{sec66} Decomposition
number $[\Delta(\lambda) : L(\mu)]$ equals zero whenever $\mu
\not \leq \lambda$. Suppose that $\mu \leq \lambda$ are
partitions in $\Lambda^+(n,r)$ such that $\lambda_1 =\mu_1$.
Define the saturated set $\Pi$ to be the set of all partitions
which are less than or equal to $\lambda$ in the dominance
ordering. In particular $\mu \in \Pi$.  We now apply the
complement construction twice. The first time we choose
$m=\lambda_1=\mu_1$. Then we obtain partitions with the $n$th
part equal to zero. We then change from $GL_n$ to $GL_{n-1}$ (see
\cite{Green}, Section 6.5) dropping the $n$th part, and apply
again the complement construction with the same $m$ as before.
Finally we change back from $GL_{n-1}$ to $GL_n$ by adding an
$n$th component zero:
\begin{eqnarray*}
[\Delta(\lambda) : L(\mu)]& =&
[\Delta(\lambda_1, \ldots ,\lambda_n) : L(\mu_1 \ldots ,\mu_n)] \\& =&
[\Delta(m-\lambda_n, \ldots ,m-\lambda_1) : L(m-\mu_n, \ldots ,m-\mu_1)] \\ &=&
[\Delta(m-\lambda_n, \ldots ,m-\lambda_2,0) : L(m-\mu_n, \ldots ,m-\mu_2,0)] \\
&=&[\Delta(m-\lambda_n, \ldots ,m-\lambda_2) : L(m-\mu_n, \ldots ,m-\mu_2)]\\
&=&[\Delta(\lambda_2, \ldots ,\lambda_n) : L(\mu_2, \ldots ,\mu_n)] \\
&=&[\Delta(\lambda_2, \ldots ,\lambda_n,0) : L(\mu_2, \ldots ,\mu_n,0)]
\end{eqnarray*}
Before we construct the isomorphism between generalised Schur algebras
underlying row removal, we first collect two remarks used in the
construction.

\subsection{Cosaturated sets and idempotents} \label{remark1} We
collect first some general facts. Let $A$ be an algebra. Given an
idempotent $e \in A$ and an ideal $J$ of $A$; by abuse of
notation we write $e$ for the image of $e$ under the natural
epimorphism $A \rightarrow A/J$. Then
\[
e(A/J)e =(eAe)/(eJe).
\]
Let $(A,\Lambda)$ be a quasi-hereditary algebra.  For a subset
$\Lambda^* \subseteq \Lambda$, we will say that $e$ is an
idempotent corresponding to $\Lambda^*$, if $e$ is the sum of
primitive orthogonal idempotents (possibly with higher multiplicities)
corresponding to the labels in $\Lambda^*$.  A subset $\Lambda^*
\subseteq \Lambda$ is called {\it cosaturated} if $\Lambda \backslash
\Lambda^*$ is saturated.

Let $e$ be the idempotent corresponding to a cosaturated set
$\Lambda^* \subseteq \Lambda$.
In this situation, $eAe$ is again quasi-hereditary with respect to
$\Lambda^*$. We call such an idempotent $e$ cosaturated.
So suppose $e$ is cosaturated in a quasi-hereditary algebra $(A,\Lambda)$
and $J$ is an ideal in the hereditary chain of $A$; this means
$A/J$ is quasi-hereditary with respect to some saturated subset
$\Lambda' \subseteq \Lambda$. Then $e$ is cosaturated in $A/J$. This
means $e(A/J)e$ is again quasi-hereditary with indexing set
$\Lambda' \cap \Lambda^*$.

\subsection{Relating $GL_n$-representations and
  $GL_{n-1}$-representations} \label{remark2} Notations and details
not given here can be found in Green \cite{Green}, in particular we
follow Section 6.5.  Given $N \geq n$. Then Green considers $I(n,r)$
as a subset of $I(N,r)$ in the natural way.  With this convention,
$S(n,r)$ can be considered as a subalgebra of $S(N,r)$. Define an
injective map from $\Lambda(n,r)$ into $\Lambda(N,r)$ by adding $N-n$
empty rows:
\[
\lambda \mapsto \lambda^*=(\lambda_1, \ldots, \lambda_n, 0
\ldots, 0).
\]
Denote the image of  $\Lambda(n,r)$ in $\Lambda(N,r)$ by
$\Lambda(n,r)^*$.
Using the notation as in Green \cite{Green}, define the idempotent
\begin{eqnarray} \label{eqn55}
   e=\sum_{\beta \in \Lambda(n,r)^*} \xi_{\beta}
\end{eqnarray}
in $S(N,r)$.  Denote by $\bar{\beta}$ the partition associated to
the composition $\beta$. Then the idempotent $\xi_{\beta}$ is
associated to $\xi_{\bar{\beta}}$. Moreover, $\xi_{\bar{\beta}}$
is a sum of primitive orthogonal idempotents corresponding to
weights greater than or equal to $\bar{\beta}$, see \cite{Green},
Section 4.7 (a).
Note that the set $\Lambda^+(n-1,r)^*$ of partitions in
$\Lambda(n-1,r)^*$ is a cosaturated subset of $\Lambda(N,r)$.
Hence $e$ is cosaturated and $eS(N,r)e$ is quasi-hereditary with
respect to $\Lambda^+(n,r)^*$ and the dominance ordering.  Moreover,
Green shows that
\begin{eqnarray}  \label{eqn33}
S(n,r) \simeq eS(N,r)e.
\end{eqnarray}
As there is an isomorphism between the indexing sets of these two
quasi-hereditary algebras, this isomorphism identifies the
quasi-hereditary structures of $eS(N,r)e$ with the quasi-hereditary
structure of $S(n,r)$.  In particular, $eL(\lambda^*)=L(\lambda)$ and
$eL(\mu)=0$ for $\mu \in \Lambda(N,r)\backslash \Lambda(n,r)^*$, and
$e\Delta(\lambda^*)=\Delta(\lambda)$.

\subsection{Row removal algebraically} We finally construct the
isomorphism between generalised Schur algebras, along the lines
of the calculation in Section \ref{sec66}. The construction is in
four steps:

(i) We first apply Theorem \ref{maintheorem}, combined with the
arguments given in Sections \ref{remark1} and \ref{remark2}. Let
$\Lambda$ be the saturated subset of $\Lambda^+(n,r)$ consisting of
all partitions with first rows of length smaller than or equal to $m$.
Then by Theorem \ref{maintheorem}:
\begin{eqnarray} \label{eqn66}
S(n,r)/I_{\Lambda} \cong S(n,nm-r)/I_{\widehat{\Lambda}}.
\end{eqnarray}
Here the indexing set $\widehat{\Lambda}$ is a saturated subset of
$\Lambda^+(n,nm-r)$, consisting of all partitions of $nm-r$ fitting
into an $n \times m$ rectangle.

(ii) The argument in Section \ref{remark2} provides an isomorphism
\begin{eqnarray} \label{eqn54}
eS(n,nm-r)e \simeq S(n-1,nm-r)
\end{eqnarray}
of quasi-hereditary algebras, where $e$ is defined like in Equation
(\ref{eqn55}): the indexing set corresponding to $e$ is
$\Lambda^+(n-1,nm-r)^*$. Let $\Pi \subseteq \Lambda^+(n-1,nm-r)$ be
such that $\Pi^*=\widehat{\Lambda} \cap \Lambda^+(n-1,nm-r)^*
\subseteq \Lambda^+(n,nm-r)$. Then $\Pi$ is a cosaturated subset of
$\Lambda^+(n-1,nm-r)$.  Hence we have the following isomorphism of
quasi-hereditary algebras:
\begin{eqnarray*}
e(S(n,r)/I_{\Lambda})e
&\cong & e(S(n,nm-r)/I_{\widehat{\Lambda}})e  \mbox{\hspace{1.65cm} by Equation (\ref{eqn66}),}\\
&\cong & (eS(n,nm-r)e)/(eI_{\widehat{\Lambda}}e)  \mbox{\hspace{1cm} by Section \ref{remark1},}\\
& \cong & S(n-1,nm-r)/(eI_{\hat{\Lambda}}e) \mbox{ \hspace{0.9cm}
by Equation (\ref{eqn54})} \\
& \cong &  S(n-1,nm-r)/ I_{\Pi}
\end{eqnarray*}
Note that in the first isomorphism, by abuse of notation, we write $e$
for the preimage of $e$ under the isomorphism in Equation
(\ref{eqn66}).
The indexing set of the left-hand side now consists precisely of the
partitions of $r$ with first row of length $m$, fitting into an $n
\times m$-rectangle. The indexing set of the right-hand side is $\Pi$.
The above isomorphism induces on the indexing sets the following
identification: given $\lambda \in \Lambda^+(n,r)$ with $\lambda_1=m$,
then
\begin{eqnarray*}
\lambda=(\lambda_1, \ldots, \lambda_n) &\rightarrow &
(m-\lambda_n, \ldots, m -\lambda_2,m -\lambda_1) \\
 &\rightarrow & (m-\lambda_n, \ldots, m -\lambda_2,0) \\
&\rightarrow & (m-\lambda_n, \ldots, m -\lambda_2).
\end{eqnarray*}

(iii) We next apply Theorem~\ref{maintheorem} with respect to $\Pi$ and
with $m$ as above.  This means we take the complement of a partition
with $nm-r$ boxes inside an $(n-1) \times m$ rectangle, and as such
the complement has $(n-1)m-(nm-r)=r-m$ boxes. Hence, by
Theorem~\ref{maintheorem}, we have
\begin{eqnarray}  \label{eqn333}
S(n-1,nm-r)/I_{\Pi} \cong S(n-1,r-m)/I_{\hat{\Pi}}.
\end{eqnarray}
On the indexing set the following identification happens:
\[
(m-\lambda_n, \ldots, m -\lambda_2) \mapsto (\lambda_2, \ldots,
\lambda_n).
\]

(iv) We finally apply again Green's isomorphism (see Section
\ref{remark2}): there exists an idempotent $f$ defined like in
Equation (\ref{eqn55}) such that
\begin{eqnarray*}
  S(n-1,r-m) \cong fS(n,r-m)f.
\end{eqnarray*}
The indexing set corresponding to $f$ is $\Lambda^+(n-1,r-m)^*$.  By
(iii), the ideal $I_{\hat{\Pi}}$ is an ideal in the hereditary chain
of $S(n-1,r-m)$. Hence
\[
 S(n-1,r-m)/I_{\hat{\Pi}} \cong (f S(n,r-m) f)/ I_{\hat{\Pi}^*},
\]
where $I_{\hat{\Pi}^*}$ is the image of $I_{\hat{\Pi}}$ under the
isomorphism in Equation (\ref{eqn333}).
In this last step, the identification on the indexing sets is:
\[
(\lambda_2, \ldots, \lambda_n) \mapsto (\lambda_2, \ldots,
\lambda_n,0).
\]
In total we constructed an isomorphism of quasi-hereditary algebras
\[
e(S(n,r)/I_{\Lambda})e \cong (f S(n,r-m) f)/ I_{\hat{\Pi}}.
\]
Here the indexing set on the left-hand side consists of partitions
$\lambda$ of $r$ with at most $n$ parts such that $\lambda_1=m$. This
guarantees that the partitions of the indexing set fit into an $n
\times m$-rectangle.  The indexing set on the right-hand side consists
of partitions of $r-m$ into $n$ parts where the $n$-th part is zero
and the first part is smaller than or equal to $m$. On the indexing
sets the above isomorphism induces the identification of
$\lambda=(\lambda_1, \ldots, \lambda_n)$ with $(\lambda_2, \ldots,
\lambda_n,0)$.

\section{Application: Equating $p$-Kostka numbers}
Fix a natural number $r$ and partitions $\lambda$ and $\mu$ of $r$.
The permutation module $M^\lambda$ over $\Sigma_r$ is the module
obtained by inducing the trivial representation from the Young
subgroup $\Sigma_\lambda$ to the symmetric group $\Sigma_r$. The Young
module $Y^\mu$ is the unique indecomposable direct summand of $M^\mu$
which contains the Specht module $S^\mu$ as a submodule. Every
permutation module $M^\lambda$ is a direct sum of Young modules
$Y^\mu$ (with multiplicities) where the indices satisfy $\mu \geq
\lambda$.  The $p$-Kostka number $(M^\lambda : Y^\mu)$ is defined to
be the multiplicity of the Young module $Y^\mu$ occuring, up to
isomorphism, in a direct sum decomposition of the permutation module
$M^\lambda$. Thus we have:
\begin{eqnarray*}
M^\lambda &= & \bigoplus_{\mu \geq \lambda} (M^\lambda : Y^\mu)  Y^\mu 
\end{eqnarray*}
If $\lambda$ and $\mu$ are partitions of $r$ with not more than $n$
parts, then we can reinterpret the $p$-Kostka number in terms of the
Schur algebra $S(n,r)$. For each partition $\lambda$, there is an
idempotent $\xi_\lambda$ in the Schur algebra, as defined in
\cite{Green,Donkin}.  The space $\xi_\lambda E^{\otimes r}$, when
viewed as right $k \Sigma_r$-module, is isomorphic to the permutation
module $M^\lambda$.  Then $\xi_\lambda \in S(n,r) = End_{k
\Sigma_r}(E^{\otimes r})$ is the identity in $End_{k
\Sigma_r}(M^\lambda)$, that is, the projection onto $M^{\lambda}$.
The identity on the indecomposable direct summand $Y^\lambda$ of
$M^\lambda$ is a primitive idempotent $e_\lambda \in S(n,r)$ contained
in $\xi_\lambda$.  We can write $\sum_{\lambda} \xi_\lambda =
\sum_\lambda \epsilon_\lambda$, where $\epsilon_\lambda$ is the sum of
all primitive idempotents equivalent to $e_\lambda$, which occur in
$\sum_{\lambda} \xi_\lambda$. Hence $\xi_{\lambda} = \sum_{\mu}
\xi_{\lambda} \epsilon_{\mu} = \sum_{\mu \geq \lambda} \xi_{\lambda}
\epsilon_{\mu}$. (Note that $\sum_{\lambda} \xi_\lambda$ is not the
unit element of $S(n,r)$, since we take the sum over partitions, not
over compositions.)
This gives us:
\begin{eqnarray*}
M^\lambda &=& \xi_\lambda E^{\otimes r} = \sum_{\mu} \xi_\lambda 
\epsilon_{\mu} E^{\otimes r}.
\end{eqnarray*}
In particular we can rewrite the $p$-Kostka number $(M^\lambda :
Y^\mu)$ as a multiplicity $[\xi_\lambda : e_{\mu}]$ of the primitive
idempotents equivalent to $e_{\mu}$ in $\xi_{\lambda}$. Note that this
is also the multiplicity of $e_\mu$ occuring, up to equivalence, in
$\xi_\lambda \cdot \epsilon_\mu$.

Appropriate sums of the idempotents $\xi_{\lambda}$ generate the
ideals in a heredity chain of the Schur algebra. Suppose $J$ is an
ideal in the heredity chain of $A$. Let \,\,$\bar{ } : A \rightarrow A/J$
be the quotient map.  Then the primitive idempotent $e_{\mu}$ is sent
either to zero, or to a primitive idempotent of $A/J$. It is sent to
zero if and only if $\xi_{\mu} \in J$. Moreover the image is non-zero,
if and only if the image of any primitive idempotent equivalent to
$e_\mu$ is non-zero, too; then the multiplicity of $e_{\mu}$ in
$\epsilon_{\mu}$ is preserved under the quotient map. Hence
\[
[\xi_\lambda : e_\mu ] = [\overline{\xi_\lambda} : \overline{e_\mu}].
\]
The isomorphism in Theorem \ref{maintheorem} preserves heredity chains
and thus sends the equivalence class of $\overline{\xi_{\lambda}}$ to
the equivalence class of $\overline{\xi_{\hat{\lambda}}}$ for $\lambda
\in \Pi$.  This proves the following corollary (which also follows
by rewriting the multiplicity $(M^\lambda : Y^\mu)$ as the dimension
of the weight space $L(\mu)^{\lambda}$ and then using the isomorphisms
$L(\mu + \omega) \simeq det \otimes L(\mu)$ and $L(\mu^{\ast}) \simeq
L(\mu)^{\ast}$, see also~\cite{Donkin}, 4.4(1)(v)):
\begin{cor}
Suppose $\lambda$ and $\mu$ are partitions in $\Lambda^+(n,r)$. Choose
a natural number $m$ such that $\lambda_1, \mu_1 \leq m$. Then 
\[
(M^\lambda : Y^\mu) = (M^{\hat{\lambda}} : Y^{\hat{\mu}}).
\]
Consequently, the column and row removal formulas for $p$-Kostka
numbers hold true.
\end{cor}
\begin{pf}
We have already shown the multiplicity formula. Note that we have
independence of $n$ in the following sense: The partitions $\lambda$ and
$\mu$ may as well be considered as elements of $\Lambda^+(N,r)$ for any
$N \geq n$ (by formally adding zeroes at the end). This does not change
the multiplicities. For instance, we have 
\[
(M^\lambda : Y^\mu) = (M^{(\lambda,0)} : Y^{(\mu,0)}).
\]
Now the row removal formula follows by applying the complement
construction twice as described in Section \ref{sec5}.

The column removal rule for $p$-Kostka numbers follows by applying the
complement construction twice as described in Section \ref{sec4}.
\end{pf}

\nocite{*} \bibliographystyle{plain}

\bigskip
{\footnotesize
\begin{tabular}{lcr}
Ming Fang & Anne Henke & Steffen Koenig \\
Dept of Mathematics \phantom{xxxxxx} & \phantom{xxx} Mathematical Institute
\phantom{xxx} & \phantom{xxxxxx} Mathematisches Institut \\
Chinese University of & University of & Universit\"at  \\
Science and Technology &Oxford & zu K\"oln \\
(USTC) & 24-29 St Giles' & Weyertal 86-90 \\
Hefei, 230026, P.R. China & Oxford OX1 3LB, UK & 50931 K\"oln, Germany \\
{\tt fming@mail.ustc.edu.cn} & {\tt henke@maths.ox.ac.uk} &
{\tt skoenig@math.uni-koeln.de} \\
{\tt  mail.ustc.edu.cn/\~{ }fming} & {\tt www.maths.ox.ac.uk/\~{ }henke}
& {\tt www.mi.uni-koeln.de/\~{ }skoenig}
\end{tabular}

}

\def\andersenpolowen{{\bf1}}
\def\donkina{{\bf2}}
\def\donkinb{{\bf3}}
\def\donkint{{\bf4}}
\def\donkinx{{\bf5}}
\def\donkinc{{\bf6}}
\def\green{{\bf7}}
\def\jantzen{{\bf8}}
\def\springer{{\bf9}}


\newcommand{\nc}{\newcommand}

\nc{\ed}{\end{document}}

\nc{\q}{\quad}
\nc{\ds}{\displaystyle}
\nc{\proj}{{\rm proj}}
\nc{\ind}{{\rm ind}}
\nc{\ch}{{\rm ch}}
\nc{\character}{{\rm char}}
\nc{\trace}{{\rm trace}}
\nc{\GL}{{\rm GL}}
\nc{\Sym}{{\rm Sym}}
\nc{\Aut}{{\rm Aut}}
\nc{\de}{\delta}
\nc{\ep}{\epsilon}
\nc{\qq}{\quad\quad}
\nc{\te}{\otimes}
\nc{\tensor}{\otimes}
\nc{\ten}{\otimes}
\nc{\op}{{\rm op}}
\nc{\rmif}{{\rm if\ }}
\nc{\aff}{{\rm aff}}
\nc{\rmotherwise}{{\rm otherwise}}
\nc{\en}{{\mathbb N}}
\nc{\eno}{{\mathbb N}_0}
\nc{\que}{{\mathbb Q}}
\nc{\real}{{\mathbb R}}
\nc{\seee}{{\mathbb C}}
\nc{\barC}{{\bar C}}
\nc{\barde}{{\bar\de}}
\nc{\barep}{{\bar\ep}}
\nc{\efp}{{\mathbb F}}
\nc{\Ker}{{\rm Ker}}
\nc{\zed}{{\mathbb Z}}
\nc{\gothg}{{\goth g}}
\nc{\gothh}{{\goth h}}
\nc{\U}{{\rm U}}
\nc{\hU}{{\hat {\rm U}}}
\nc{\hatU}{\hU}
\nc{\dash}{---}
\nc{\bs}{\bigskip}
\nc{\ssl}{{\rm sl}}
\nc{\Ext}{{\rm Ext}}
\nc{\alg}{{\rm alg}}
\nc{\Prim}{{\rm Prim}}
\nc{\Fin}{{\rm Fin}}
\nc{\Inv}{{\rm Inv}}
\nc{\triv}{{\rm triv}}
\nc{\Ad}{{\rm Ad}}
\nc{\Cent}{{\rm Cent}}
\nc{\Loc}{{\rm Loc}}
\nc{\Ann}{{\rm Ann}}
\nc{\Mod}{{\rm Mod}}
\nc{\mathcalB}{{\mathcal B}}
\nc{\mathcalC}{{\mathcal C}}
\nc{\F}{{\mathcal F}}
\nc{\Rep}{{\rm Rep}}
\nc\cards{\quad THIS MUST GO\quad}

\nc{\ba}{\begin{align*}}

\nc{\Comod}{{\rm Comod}}
\nc{\la}{\lambda}
\nc{\La}{\Lambda}
\def\mathcalL{{\mathcal L}}
\def\mathcalR{{\mathcal R}}
\nc{\cf}{{\rm cf}}
\nc{\Gp}{{\rm Gp}}
\nc{\union}{\cup}
\nc{\iso}{\cong}
\nc{\Pf}{{\bf Proof.\rm\ }}
\nc{\cl{\centerline}}

\parindent=0pt

\newpage

{\Large{\centerline{\bf Appendix by Stephen Donkin}}}

\thispagestyle{empty}

\bs\bs

\rm\q We give here an alternative approach which gives substantial
generalizations of the main result, Theorem \ref{maintheorem} above.

\q We first consider the classical case $q=1$.  The Schur algebra
$S(n,r)$, as defined by Green, [\green], is the dual algebra of the
subcoalgebra $A(n,r)$ of the coordinate algebra of a the general
linear group $\GL_n$ and a similar construction of "generalized Schur
algebras" is given in [\donkina], [\donkinb] in the context of
arbitrary reductive groups. It therefore seems natural to proceed
directly in terms of coalgebras. The dual of an isomorphism of
coalgebras is an isomorphism of the dual algebras and so we, for the
most part, restrict ourselves to producing coalgebra isomorphisms.

\q We fix a connected reductive group $G$ defined over a field $k$
which we assume algebraically closed (but see Remark 4 below).  We
regard $k[G]$ as a Hopf algebra over $k$ with structure maps
$\de:k[G]\to k[G]\otimes k[G]$, $\ep:k[G]\to k$, $\sigma:k[G]\to k[G]$
(as in [\donkinc; 0.7] for example)

\q For an algebraic group $H$ over $k$ we write $X(H)$ for the additive 
abelian group  of algebraic group  homomorphisms from $H$ to the 
multiplicative group of $k$. We assume the usual notation for reductive 
groups and their representations, as in for example [\donkina],[\donkinb].  
In  particular, we have, inside the weight lattice $X=X(T)$,  the set of 
dominant weights $X^+(T)$.  For each $\lambda\in X^+(T)$ we have the 
irreducible module $L(\lambda)$ of highest weight $\lambda$, and the module 
$\nabla(\lambda)$ induced from a one dimensional $B$-module on which $T$ 
acts via $\lambda$. Here $B$ is a Borel subgroup whose set of roots forms a 
system of negative roots in the root system of $G$, with respect to $T$.   
We write $\Phi^+$ (resp. $\Phi^-$) for the set of positive (resp. negative) 
roots and $\Phi$ for $\Phi^+\cup \Phi^-$, the set of roots. 
There is a natural partial order on $X(T)$: we write $\lambda\leq \mu$ if 
$\mu-\lambda$ is a sum of positive roots. For a $G$-module $V$ we write 
$\cf(V)$ for the coefficient space of $V$, i.e. the subspace of $k[G]$ 
spanned by the coefficient functions $f_{ij}$, $i\in I$, where $v_i$, 
$i\in I$, is a basis of $V$ and the $f_{ij}$ are defined by the equations
$$gv_i=\sum_{j\in I}f_{ji}(g)v_j$$
for all $i\in I$,  $g\in G$. Recall that for $G$-modules $V_1,V_2$ we have
$\cf(V_1\otimes V_2)=\cf(V_1)  \cdot \cf(V_2)$.

\q We recall the construction of the generalized Schur algebras.  Let $\Pi$ 
be a subset of $X^+(T)$. A (rational) $G$-module $V$ is said to belong to 
$\Pi$ if every composition factor of $V$ belongs to 
$\{L(\lambda)\vert\lambda\in\Pi\}$. 
Among all submodules belonging to $\Pi$ of an arbitrary $G$-module $V$
there is a unique maximal one, denoted $O_\Pi(V)$.  In this way we get a 
left exact functor $O_\Pi$ from the category of rational $G$-modules to 
itself. The set $\Pi$ is said to be saturated if it has the property that 
whenever $\lambda\in\Pi$ and $\mu$ is a dominant weight less such that 
$\mu\leq \lambda$ then $\mu\in\Pi$. Regarding $k[G]$ as a rational left 
$G$-module we have the submodule $A(\Pi)=O_\Pi(k[G])$, and indeed $A(\Pi)$ 
is a subcoalgebra of $k[G]$. In fact $A(\Pi)$ is the sum of all spaces 
$\cf(V)$, as $V$ ranges over all $G$-modules belonging to $\Pi$.

\q It is shown in [\donkina] that if
$\Pi$ is finite and saturated then $A(\Pi)$ is  finite
dimensional. However, saturation is not important for what follows so we
now allow $\Pi$ to be any finite subset or $X^+(T)$. As  $A(\Pi)\subset
A(\Gamma)$, whenever $\Pi\subseteq \Gamma\subseteq X^+(T)$ we have that
$A(\Pi)$ is finite dimensional since we may take for $\Gamma$ the smallest
saturated subset containing $\Pi$, and this is finite. The generalized
Schur algebra $S(\Pi)$ is defined as the dual algebra of the coalgebra  
$A(\Pi)$. 

\q We shall describe two natural situations in which one gets an
isomorphism of Schur coalgebras $A(\Gamma)\to A(\Pi)$ (and hence of Schur 
algebras $S(\Pi)\to S(\Gamma)$). First suppose that $L$ is a one
dimensional $G$-module and that $\mu\in X(T)$ is the representation of $T$ 
afforded by $L$. Then $L=L(\mu)$. We let $d_\mu$ be the element of $k[G]$ 
such that $gx=d_\mu(g)x$, for all $g\in G$, $x\in L$.  Notice that the dual 
module $L(\mu)^*$ is isomorphic to $L(-\mu)$. It follows that $d_\mu 
d_{-\mu}=1$, in particular $d_\mu\in k[G]$ is a unit.

\q Let $\Pi$ be a subset of $X^+(T)$. We have $L\otimes L(\lambda)\iso 
L(\lambda+\mu)$, for $\lambda\in X^+(T)$. Hence a $G$-module $V$ has 
composition factor $L(\lambda)$ if and only if $L\otimes V$ has composition 
factor $L(\lambda+\mu)$.  If $V$ is a $G$-module belonging to $\Pi$ then 
$L\otimes V$ belongs to $\mu+\Pi=\{\lambda+\mu\vert\lambda\in \Pi\}$ and 
hence  $d_\mu\cf(V)=\cf(L\otimes V)\subseteq A(\mu+\Pi)$. 
Thus we get $d_\mu A(\Pi)\subseteq A(\mu+\Pi)$. But, by the same token, we 
have $d_{-\mu}A(\mu+\Pi)\subset A(\Pi)$ and hence $d_\mu
A(\Pi)=A(\mu+\Pi)$. Since $d_\mu(gh)=d_\mu(g)d_\mu(h)$, for all $g,h\in G$, 
we have $\de(d_\mu)=d_\mu\otimes d_\mu$ and we get our first source of 
isomorphisms of Schur coalgebras.
\bs

\bf Principle 1\q\sl Suppose that $\mu\in X(T)$ is such that $L(\mu)$ is
one dimensional and $d_\mu\in k[G]$ is given by $gx=d_\mu(g)x$ for all
$g\in G$, $x\in L(\mu)$. Then the map  $A(\Pi)\to A(\mu+\Pi)$, 
given by multiplication by $d_\mu$, is a coalgebra isomorphism.

\thispagestyle{empty}

\bs\rm

\bf Remark 1\q\rm We thus get an isomorphism of Schur coalgebras for each
element of $X(G/G')$, where $G'$ is the derived group of $G$.  As we have 
$G=TG'$ we may identify $X(G/G')$ with $X(T/T')$, where $T'=G'\cap T$. In 
other words we have specifically an isomorphism $A(\Pi)\to A(\mu+\Pi)$ for 
each  $\mu\in X(T)$ which has $T'$ in its kernel.

\bs

\bf Remark 2\q\rm  The arguments above are formal and hold  quite
generally.  Let $(H,\de,\ep)$ be a Hopf  algebra. Let $\Gp(H)$ be the set
of group-like elements, i.e. an elements $g\in H$ such $d(g)=g\otimes g$ 
and  $\ep(g)=1$. Then $\Gp(H)$ is a subgroup of the group of units of $H$.
Let $\{L_\lambda\vert\lambda\in X\}$ be a complete set of pairwise 
non-isomorphic right $H$-comodules. 
There is a natural action of $\Gp(H)$ on $X$ defined as follows.  If $g\in 
\Gp(H)$ then $kg$ is a one dimensional subcomodule of $H$ and, for 
$\lambda\in X$, the comodule $kg\otimes L_\lambda$ is simple. We define 
$g\lambda\in X$ by the condition $kg\otimes L_\lambda\iso L_{g\lambda}$.

\q For a subset $\Pi$ of $X$, we may define a functor $O_\Pi$, form
the category of right $H$-comodules to itself. Then viewing $H$ itself
as a right $H$-comodule (with structure map $\de:H\to H\otimes H$) we
may form the subcoalgebra $A(\Pi)=O_\Pi(H)$. The above argument shows
that multiplication by $g$ determines an isomorphism of coalgebras
from $A(\Pi)$ to $A(g\Pi)$.

\bs

\q Our second source of isomorphisms of Schur algebras is from
automorphisms of $G$. Let $\phi:G\to G$ be an automorphism of our reductive 
group $G$. For a $G$-module $V$ affording the representation $\pi:G\to 
\GL(V)$ we form the module $V^\phi$ with the same underlying $k$-space $V$ 
on which $G$ acts according to the representation $\pi\circ\phi$. We note 
that $(V^\phi)^\psi\iso V^{\phi\circ\psi}$, for automorphisms $\phi,\psi$, 
and that $V\iso V^\phi$ if $\phi$ is inner. Moreover, we have (from the 
definitions)  $\cf(V^\phi)=\phi^\sharp(\cf(V))$, for a $G$-module $V$, 
where  $\phi^\sharp:k[G]\to k[G]$ is the comorphism of $\phi$. 
Hence $\phi^\sharp(A(\Pi))=A(\Pi)$ if $\phi$ is an inner automorphism. Let 
$\phi$ be a general automorphism. Then  by the conjugacy theorems for Borel 
subgroups and maximal tori, there exists an element $g\in G$ such that 
$\phi(B)=B^g$ and $\phi(T)=T^g$. Thus we may write $\phi=\gamma\circ\psi$, 
where $\psi$ is an automorphism stabilizing $B$ and $T$ and $\gamma$ is an 
inner automorphism.
Now $\psi$ restricts to an isomorphism $\psi_0:T\to T$ which induces an 
isomophsm $f : X(T)\to X(T)$ (the restriction of the comorphism 
$\psi_0^\sharp:k[T]\to k[T]$).  Moreover, since $\psi$ preserves $B$, the 
isomorphism $f$ induces a  bijection on negative roots, and hence also on
 positive roots.  The map $f$ is an automorphism of the root  datum 
$(X,\Phi,\check X, \check \Phi)$ and conversely, for any automorphism of
the root datum $f'$ there is an automorphism $\psi':G\to G$ stabilizing $B$ 
and $T$ and inducing $f$ on $X(T)$ (see [\springer; 11.4.3 Theorem]).    

\q Now let $\Pi$ be a finite subset of $X^+(T)$ and let $V$ be a module 
belonging to $\Pi$. Then the weights of $V^\phi$ are the weights of the
form $f(\mu)$, for $\mu$ a weight of $V$. Hence $V^\phi$ belongs to
$f(\Pi)$ and $\cf(V^\phi)=\phi^\sharp(cf(V))\leq A(f(\Pi))$. Thus we get 
$\phi^\sharp(A(\Pi))\leq A(f(\Pi))$ and hence (applying the same principle 
with the inverse of $\psi$ in place of $\psi$) we have 
$\phi^\sharp(A(\Pi))=A(f(\Pi))$.  
\bs

\bf Principle 2\q\sl Let $\Pi$ be a finite subset of $X^+(T)$. Let $\phi$ be an isomorphism of $G$ inducing the isomorphism $f$ on the root datum as above. Then the restriction of $\phi^\sharp$ is an isomorphism from $A(\Pi)$ to $A(f(\Pi))$. Conversely, for any automorphism of the root datum $f$ there is an isomorphism $A(\Pi)\to A(f(\Pi)$, namely the restriction of $\phi^\sharp$ where $\phi$ is an automorphism of $G$ stabilizing $B$ and $T$ and inducing $f$ on $X(T)$.

\thispagestyle{empty}

\bs

\bf Remark 3\q\rm  We may define $f:X\to X$ by $f(\mu)=-w_0\mu$. For 
$\lambda$ a dominant weight we set $\lambda^*=-w_0\lambda=f(\lambda)$,
where $w_0$ is the longest element of the Weyl group $W=N_G(T)/T$,  and for 
a subset $\Pi$ of $X^+(T)$ set $\Pi=f(\Pi)=\{-w_0\lambda\vert\lambda\in
\Pi\}$.  Thus, for each finite subset $\Pi$ of $X^+(T)$,  we get a natural 
isomorphism $A(\Pi)\to A(\Pi^*)$.

\bs

\bf Remark 4\q \rm Suppose  $G$ is defined and split over an arbitrary
field $k$.  The  first principle is still valid and the proof goes through 
without change to give isomorphisms of the $k$-coalgebras $A(\Pi)$ and 
$A(\mu+\Pi)$ (defined as subcoagebras of $k[G]$).  Moreover, for each 
automorphism $f$ of the root datum there is a $k$-isomorphism $\phi:G\to G$ 
inducing $f$ (see e.g.  [\jantzen;II, 1.15 Proposition]) so that the second 
principle is valid also in this case and we have an isomorphism of the 
coalgebras $A(\Pi)$ and $A(f(\Pi))$ (defined as subcoalgebras of $k[G]$).

\bs

\bf Remark 5 \q\rm We shall not attempt to give general versions of 
principles of quantum groups of Principles 1 and 2, though it is clear that 
it would be possible to do so within the framework of, for example, 
[\andersenpolowen]. 

\bs

\q We now turn our attention to the quantum general linear group $G$, of 
degree $n$,over an arbitrary field $k$, with parameter $0\neq q\in k$,  as 
in [\donkinc]. We have the set of dominant weights $X^+(n)$ and the set of 
polynomial dominant weights $\Lambda^+(n)$ and, for $r\geq 0$, the set 
$P^+(n,r)$ of polynomial dominant weights $\lambda=(\lambda_1,\ldots,
\lambda_n)$ such that $\lambda_1+\cdots+\lambda_n=r$. For $\lambda=
(\lambda_1,\ldots,\lambda_n)\in X^+(n)$ we put $\lambda^*=(-\lambda_n,
\ldots,-\lambda_1)$ (i.e. $-w_0\lambda$, where $w_0$ is the longest element 
of the Weyl group). For $\Pi\subset X^+(T)$, we put $\Pi^*=
\{\lambda^*\vert\lambda\in\Pi\}$. For $\Pi$ finite (and not necessarily 
saturated) we have the  Schur coalgebra $A(\Pi)$ and its dual algebra 
$S(\Pi)$ defined as above (see also   [\donkinx]). We shall produce an 
isomorphism $A(\Pi)\to A(\Pi^*)$, for $\Pi$ finite. (Note that in the case 
$q=1$ we could directly  invoke Principle 2.)  Note that multiplication by 
the determinant $d\in k[G]$ gives an isomorphism $A(\Pi)\to dA(\Pi)=
A(\omega+\Pi)$ as in Principle 1, where $\omega=(1,1,\ldots,1)$.

\q We now show that $A(\Pi)$ is isomorphic to $A(\Pi^*)$. Let
$\sigma:k[G]\to k[G]$ be the antipode. Then the relationship between the
coefficient spaces of a finite dimensional left $G$-module $V$ and the dual 
left $G$-module $V^*$ is $\cf(V^*)=\sigma(\cf(V))$. Now if $V$ is such that 
$\cf(V)=A(\Pi)$ (e.g. $V$ is $A(\Pi)$ itself as a left $G$-module) then we 
get $\sigma(A(\Pi))=\sigma(\cf(V))=\cf(V^*)\leq A(\Pi^*)$. Thus we get 
$\sigma^2(A(\Pi))\leq \sigma(A(\Pi^*)\leq A(\Pi)$ and now by dimensions
(and the fact that $\sigma$ is injective) we have
$\sigma(A(\Pi))=A(\Pi^*)$. Moreover $\sigma : k[G]\to k[G]$ is an 
antimorphism of  coalgebras  so we get that $\sigma$ induces an isomorphism 
$S(\Pi^*)\to S(\Pi)^\op$ (where ${}^\op$ indicates the opposite algebra). 
[Note that this is a general argument for Hopf algebras and subcoalgebras 
defined by restricting composition factors.]  It remains to prove that 
$S(\Pi)$ is isomorphic to $S(\Pi)^\op$ for any finite saturated subset
$\Pi$. We can replace $S(\Pi)$ by the isomorphic algebra $S(\Pi+m\omega)$ 
and choosing $m$ large we can assume that $\Pi$ consists of polynomial 
weights. Then $A(\Pi)=\oplus_{r\geq 0}A(\Pi(r))$, where $\Pi(r)=\Pi\cap 
P^+(n,r)$, so we may assume $\Pi$ homogeneous of some degree.  But there is 
an anti-automorphism $J$ of $S(n,r)$ (see p82 of [5]) which fixes each 
idempotent $\xi_\alpha$.  By means of $J$ one generalizes to general $q$
the contravariant dual $M^0$ of a finite dimensional $S(n,r)$-module $M$
(as discussed by Green, [\green],  in the case $q=1$) - the action is 
$(a\alpha)(m)=\alpha(J(a)m)$, for $\alpha\in M^0=\Hom_k(M,k)$, $m\in M$, 
$\alpha\in M^0$.
Moreover, $M$ and $M^0$ have the same character hence the same composition 
factors. Now if $a$ annihilates $M$ then $J(a)$ annihilates $M^0$, so 
$J(I_\Pi)\leq I_\Pi$ and $J$ induces an isomorphism  $S(\Pi)\to S(\Pi)^\op$. 

 \q This brings us to our third principle, which is a combination of 
principles 1 and 2 in the situation of quantum general linear groups. 

\bs

\thispagestyle{empty}

\bf Principle 3\q\sl For any finite subset $\Pi$ of $X^+$ and $m\in\zed$ we have 
$A(\Pi)\iso A(m\omega+\Pi^*)$ and hence $S(\Pi)\iso S(m\omega+\Pi^*)$.

\bs
\rm

\bf Remark 6\q\rm We fix $r$ and $m$ and let $\Pi$ be any subset of 
$P^+(n,r)$ consisting of partitions
$\lambda=(\lambda_1,\ldots,\lambda_n)\in P^+(n,r)$ with $\lambda_1\leq m$. 
Then $\hat\Pi$ (as in Theorem \ref{maintheorem}) is $m\omega+\Pi^*$. Hence 
we get $S(\Pi)\iso S(\hat\Pi)$.
If $\Gamma\subset \Sigma \subset X^+$ then we get (from the definitions) 
that $A(\Gamma)\subset A(\Sigma)$ and the kernel $I_{\Gamma,\Sigma}$, say,  
of the  surjective algebra homomorphism $S(\Sigma)\to S(\Gamma)$ consists
of all $x\in S(\Sigma)$ which annihilate all modules belonging to
$\Gamma$. Since $\Pi\subset P^+(n,r)$ and $\hat\Pi\subset P^+(n,nm-r)$ we 
get $S(n,r)/I_\Pi\iso S(n,nm-r)/I_{\hat\Pi}$, where 
$I_{\Pi}=I_{\Pi,P^+(n,r)}$ and $I_{\hat\Pi}=I_{\hat\Pi,P^+(n,nm-r)}$. This 
gives the main result of this paper, Theorem \ref{maintheorem}, but without 
the restriction of saturation.

\bs

\bf Remark 7 \q\rm  Finally, we remark that the isomorphism above are 
defined integrally for $\Pi$ saturated. If $G$  is a general linear group
or Chevalley group, then $G$ is defined over $\zed$ 
in the usual way. This amounts to giving a   suitable $\zed$-form of the 
coordinate algebra of the complex group and gives rise to integral Schur 
coalgebras $A(\Pi)_\zed$
 (see  [\donkina; Section 4] and [\donkint])  and the Schur algbra over an 
arbitrary field is obtained by base change from the integral one. 
 Then the above arguments are valid over $\zed$ and give isomorphisms 
$S_\zed(\Pi)\to S_\zed(\Pi+\lambda)$, $S(\Pi)\to S(\Pi^*)$ etc. So one 
obtains 
$S_\zed(n,r)/I_{\Pi,\zed}\to \ldots$  (where $I_{\Pi,\zed}=
I_{\Pi,\seee}\cap S_\zed(n,r)$), and in the quantum case  
$S(n,r)_{\zed[t,t^{-1}]}/I_{\zed[t,t^{-1}],\Pi}\to...$ which specializes to 
Theorem \ref{maintheorem} over a field by base change.  
\bs

\thispagestyle{empty}

\bigskip

\footnotesize 
Stephen Donkin, Department of Mathematics,  University of
York,  Heslington,  York YO10 5DD.\\
Email: {\tt sd510@york.ac.uk}

\end{document}